\documentclass[12pt]{article}
\usepackage[mathscr]{euscript}
\usepackage{amsmath,amssymb}
\usepackage{verbatim}
\usepackage{dsfont,bm}
\usepackage{color}
\usepackage{graphicx}
\usepackage{amsthm}
\usepackage{enumerate}
\usepackage{esint} %\fint

\usepackage{xcolor}
\definecolor{darkgreen}{rgb}{0,0.75,0}
\definecolor{darkred}{rgb}{0.75,0,0}
\definecolor{darkmagenta}{rgb}{0.5,0,0.5}
\usepackage[colorlinks,citecolor=darkgreen,linkcolor=darkred,urlcolor=darkmagenta]{hyperref}
\newtheorem{theorem}{Theorem}[section]

\newtheorem{cor}[theorem]{Corollary}
\newtheorem{lemma}[theorem]{Lemma}
\newtheorem{lem}[theorem]{Lemma}
\newtheorem{proposition}[theorem]{Proposition}
\newtheorem{prop}[theorem]{Proposition}

\theoremstyle{definition}
\newtheorem{definition}[theorem]{Definition}

\newtheorem{assumption}[theorem]{Assumption}
\newtheorem{remark}[theorem]{Remark}
\newtheorem{example}[theorem]{Example}

\numberwithin{equation}{section}

\def\be{\begin{equation}}
	\def\ee{\end{equation}}
\def\bes{\begin{equation*}}
	\def\ees{\end{equation*}}

\newcommand{\Mod}[0]{\operatorname{Mod}}

\newcommand{\arxiv}[1]{{\tt \href{http://arxiv.org/abs/#1}{arXiv:#1}}}
\newcommand{\set}[1]{\left\{ #1 \right\}}

\newcommand{\abs}[1]{{\left\vert\kern-0.25ex #1
		\kern-0.25ex\right\vert}}
\def\sD {{\mathcal D}}  
  
\def\sJ {{\mathcal J}}

\def\sS {{\mathcal S}}

 \def\bN {{\mathbb N}}

 \def\bZ {{\mathbb Z}}

\def\ignore#1{}

\def\lam {\lambda}  
\def\eps{\varepsilon}

\def\Gam{\Gamma} \def\gam{\gamma}

\def\to {\rightarrow}

\def\q{\quad} 
\def\dint{\int\kern-.6em\int}

\newcommand\restr[2]{{% we make the whole thing an ordinary symbol
		\left.\kern-\nulldelimiterspace % automatically resize the bar with \right
		#1 % the function
		\vphantom{\big|} % pretend it's a little taller at normal size
		\right|_{#2} % this is the delimiter
}} %restriction of a function
\def\diam{{\mathop{{\rm diam }}}}

\def\Mod{{\mathop{{\rm Mod}}}}
\newcommand{\on}[1]{\operatorname{ #1}}

\def\wt{\widetilde}
\def\wh{\widehat}

\def\be{\begin{equation}}
	\def\ee{\end{equation}}
\def\bes{\begin{equation*}}
	\def\ees{\end{equation*}}
\def\ba{\begin{align}}
	\def\ea{\end{align}}
\def\xxea{\end{align}}
\def\bas{\begin{align*}}
\def\eas{\end{align*}}

\def\qed{{\hfill $\square$ \bigskip}}

\definecolor{dgreen}{rgb}{0, 0.6, 0.1}
\definecolor{dblue}{rgb}{0, 0.0, 0.6}
\definecolor{vdblue}{rgb}{0,.08, 0.45}
\definecolor{dred}{rgb}{0.7, 0.0, 0.0}
\definecolor{vdblue}{rgb}{0,.08, 0.45}

\definecolor{purple}{rgb}{0.6, 0.0, 0.6}
\definecolor{mytext}{rgb}{0.1, 0.1, 0.1}

\def\Adm{\operatorname{Adm}}

\begin{document}

\font\titlefont=cmbx14 scaled\magstep1
\title{Conformal Assouad dimension as the critical exponent for combinatorial modulus}    
\author{Mathav Murugan\footnote{Research partially supported by NSERC and the Canada research chairs program.} 
}
\maketitle
\vspace{-0.5cm}
%\ver
\begin{abstract}
	The conformal Assouad dimension is the infimum of all possible values of  Assouad dimension after a quasisymmetric change of metric.
	We show that the conformal Assouad dimension equals a critical exponent associated to the combinatorial modulus for any compact doubling metric space. This generalizes a similar result obtained by Carrasco Piaggio for the Ahlfors regular conformal dimension to a larger family of spaces. We also show that the value of conformal Assouad dimension is unaffected if we replace quasisymmetry with power quasisymmetry in its definition. 
	\vskip.2cm
	%\noindent {\it Keywords:} 
	\
\end{abstract}
\section{Introduction}
The \emph{Assouad dimension} of a metric space $(X,d)$ is defined as
\[
\dim_{\on{A}}(X,d)= \inf \Biggl\{ \beta >0  \Biggm|
\begin{minipage}{260 pt} there exists $C>0$ such that $N_{ r} (B(x,R)) \le C \left(\frac{R}{r} \right)^\beta$ for any $x \in X, 0<r<R$
\end{minipage}
 \Biggr\},
\]
where $N_{r}(A)$ denotes the minimum number of balls of radii $r$ required to cover $A\subset X$.
 Equivalently, Assouad dimension is the infimum of all numbers $\beta>0$ such that there exists $C>0$ so that every ball of radius $r$ has at most $C \eps^{-\beta}$ distinct points whose mutual distance is at least $\eps r$ \cite[Exercise 10.17]{Hei}.  We refer to the recent book by Fraser \cite{Fra} for a comprehensive background.
 
 We recall the definition of the conformal gauge. This terminology is motivated from the understanding that quasisymmetric maps are an analogue of conformal maps in the context of metric spaces. 
 \begin{definition}[Conformal gauge] \label{d:cgauge}
 	Let $(X,d)$ be a metric space and $\theta$ be another metric on $X$. We say that $d$ is \emph{quasisymmetric} to $\theta$, if there exists a homeomorphism $\eta:[0,\infty) \to [0,\infty)$ such that
 	\[
 	\frac{\theta(x,a)}{\theta(x,b)} \le \eta\left(\frac{d(x,a)}{d(x,b)}\right)   \q \mbox{for all triples of points $x,a,b \in X$, $x \neq b$.}
 	\]
 	We say that $d$ is \emph{power quasisymmetric} to $\theta$ if the homeomorphism $\eta$ above can be chosen so that $ \eta(t)= C (t^\alpha \vee t^{1/\alpha})$ for all $t >0$, where $C,\alpha \in [1,\infty)$.
 	The \emph{conformal gauge} of a metric space $(X,d)$ is defined as
 	\begin{equation} \label{e:cgauge-dfn}
 		\sJ(X,d):= \{ \theta: X \times X \to [0,\infty) \mid \mbox{$\theta$ is a metric on $X$, $d$ is quasisymmetric to $\theta$} \}.
 	\end{equation}
 	We define the \emph{power quasisymmetric conformal gauge} of $(X,d)$ as 
 	\begin{equation} \label{e:pcgauge-dfn}
 		\sJ_p(X,d):= \{ \theta \in \sJ(X,d) \mid \mbox{$\theta$ is  power quasisymmetric to $d$} \}.
 	\end{equation}
The \emph{conformal Assouad dimension} of $(X,d)$ is defined as
 \be  \label{e:defdca}
 \dim_{\on{CA}}(X,d)= \inf \{ \dim_{\on{A}}(X,\theta): \theta \in \sJ(X,d)\},
 \ee 
 where $ \dim_{\on{A}}(X,\theta)$ denotes the Assouad dimension of $(X,\theta)$.
\end{definition}
 
 As our main result relates conformal Assouad dimension with combinatorial modulus, we recall the notion of combinatorial modulus and a critical exponent associated to it.
 The combinatorial $p$-modulus of a family of curves $\Gam$ in a graph $G=(V,E)$ is defined as 
 \[
 \Mod_p(\Gam,G)= \inf \Biggl \{ \sum_{v \in V} \rho(v)^p \mid  \rho:V \to [0,\infty), \sum_{v \in \gam} \rho(v) \ge 1 \mbox{ for all $\gam \in \Gam$}\Biggr \}.
 \]
 Fix parameters $a, \lam, L >1$. We choose a sequence $X_k , k \ge 0$  such that $X_k$ is a maximal $a^{-k}$-separated subset of $(X,d)$ and $X_k \subset X_{k+1}$ for all $k \ge 0$. For each $k$, we define a graph $G_k$ whose vertex set is $X_k$ and there is an edge between two distinct vertices $x,y \in X_k$ if and only if $B(x, \lam a^{-k}) \cap B(y, \lam a^{-k}) \neq \emptyset$. We think of $G_k$ as a sequence of combinatorial approximations of $(X,d)$ at scale $a^{-k}$.
 We define \[ M_{p,k}(L)=  \sup \{ \Mod_{p} (\Gam_{k,L}(x), G_{k+n}) \mid x \in X_n, n \ge 0 \} \mbox{ and } M_p(L) = \liminf_{k \to \infty} M_{p,k}(L),
 \]
 where $\Gam_{k,L}(x)$ is the family of paths in $G_{n+k}$ from $B(x,a^{-n})$ to $B(x,La^{-n})^c$ (see \textsection \ref{s:modulus} for a detailed definition).
 The critical exponent corresponding to combinatorial modulus is defined as
 \[
\on{CE}(X,d) = \inf \{p>0 \mid M_{p}(L) =0 \}.
 \]
 It is not difficult to show that $\on{CE}(X,d) $ is well-defined in the sense that $\on{CE}(X,d)$ does not depend on the precise choices of $a,L,\lam \in (1,\infty)$ and also on the choices of $X_k$ (see Proposition \ref{p:wd}). Since it only depends on the metric space $(X,d)$, our notation $\on{CE}(X,d)$ is justified.

Our main result is the following theorem.
\begin{theorem} \label{t:main}
	Let $(X,d)$ be a compact metric space such that $\dim_{\on{A}}(X,d)<\infty$. Then
	\[
	\dim_{\on{CA}}(X,d)=\on{CE}(X,d)= \inf \{\dim_{\on{A}}(X,\theta): \theta \in \sJ_p(X,d)\}.
	\]
\end{theorem}

 A similar result was obtained by Carrasco \cite[Theorem 1.3]{Car} for the Ahlfors regular conformal dimension and independently in an unpublished work of Keith and Kleiner. These works rely crucially on ideas of Keith and Laakso who first related conformal Assouad dimension to combinatorial modulus \cite{KL}. 
  To state Carrasco's result, we recall the definition of Ahlfors regular conformal dimension and related notions.
 A Borel measure $\mu$ on $(X,d)$ is said to be \emph{$p$-Ahlfors regular} if there exists $C \ge 1$ such that 
 \[
 C^{-1} r^p \le \mu(B(x,r)) \le C r^p \quad \mbox{for all $x \in X, 0< r \le \diam(X,d)$.}
 \]
 Note that if such a $p$-Ahlfors regular $\mu$ exists, then the $p$-Hausdorff measure is also $p$-Ahlfors regular. The Ahlfors regular conformal dimension is defined as
 \be  \label{e:defdarc}
 \dim_{ \on{ARC}}(X,d)= \inf \{ p>0: \theta \in \sJ(X,d) \mbox{ and $\mu$ is a $p$-Ahlfors regular measure on $(X,\theta)$} \}.
 \ee
 In \eqref{e:defdarc} and \eqref{e:defdca} we adopt the convention that $\inf \emptyset = \infty$. Ahlfors regular conformal dimension is a well-studied notion in complex dynamics and hyperbolic groups; see for example \cite{BK05, BM, CM, HP09, PT, Par}. These notions of conformal dimensions are  variants of the one introduced by Pansu \cite{Pan89} and we refer the reader to \cite{MT10} for more background and applications.

To compare our results with earlier ones, we recall the notion of doubling and uniformly perfect metric spaces. 
  A metric space is said to be \emph{doubling}, if there exists $N \in \bN$ such that every ball of radius $r$ can be covered by at most $N$ balls of radii $r/2$. It is easy to see that the $\dim_{\on{A}}(X,d) <\infty$ if and only if $(X,d)$ is doubling. 
 A metric space $(X,d)$ is said to be \emph{uniformly perfect} if there exists $C>1$ such that whenever $B(x,r) \neq X$, we have $B(x,r) \setminus B(x,r/C) \neq \emptyset$. Carrasco's theorem \cite[Theorem 1.3]{Car} states that for any compact, doubling, uniformly perfect metric space, the Ahlfors regular conformal dimension is given by
 \[
 \dim_{ \on{ARC}}(X,d)= \on{CE}(X,d).
 \]
 The following lemma characterizes the class of metric spaces for which $\dim_{\on{CA}}(X,d)$ and $\dim_{\on{ARC}}(X,d)$ are finite.
\begin{lem} \label{l:basic}
	Let $(X,d)$ be a compact metric space. Then 
	\begin{enumerate}[(a)]
		\item  $\dim_{\on{CA}}(X,d)$ is finite if and only if $(X,d)$ is doubling.
		\item  $\dim_{\on{ARC}}(X,d)$ is finite if and only if $(X,d)$ is doubling and uniformly perfect. Moreover,
		 If $(X,d)$ is doubling and uniformly perfect, then  $\dim_{\on{ARC}}(X,d)=\dim_{\on{CA}}(X,d)$ \cite[Proposition 2.2.6]{MT10}.
	\end{enumerate}
\end{lem}
 By Lemma \ref{l:basic},   our result in Theorem \ref{t:main}  generalizes Carrasco's  theorem \cite[Theorem 1.3]{Car} to  doubling metric spaces that are not necessarily uniformly perfect.  
We refer to \cite{Kig,Sha} for expositions to Carrasco's work.

One motivation for this work is that conformal Assouad dimension is better behaved than Ahlfors regular conformal dimension.
The above lemma shows that conformal Assouad dimension is a 
meaningful quasisymmetry invariant for a larger class of metric spaces.
If $(X,d)$ is a compact metric space  and $Y \subset X$, then it is easy to see that
\[
\dim_{\on{CA}}(Y,d) \le \dim_{\on{CA}}(X,d).
\]
The above inequality is not always true for Ahlfors regular conformal dimension because a subset  of uniformly perfect metric space is not necessarily uniformly perfect. Nevertheless, if $(X,d)$ is a compact, doubling, uniformly perfect metric space and $Y \subset X$ is also compact, doubling and uniformly perfect, then
\be  \label{e:subset}
\dim_{\on{ARC}}(Y,d) \le \dim_{\on{ARC}}(X,d).
\ee
One way to show \eqref{e:subset} is to use Lemma \ref{l:basic}(b) and the analogous inequality for conformal Assouad dimension. Another more involved approach would be to use \cite[Theorem 1.3]{Car} and   careful choices of hyperbolic fillings for $X$ and $Y$. The direct approach of restricing an Ahlfors regular metric in $\sJ(X,d)$ to $Y$ does not work because the restriction of an Ahlfors regular metric  on a subset need not be Ahlfors regular.   To summarize, conformal Assouad dimension is better behaved because Ahlfors regularity and uniform perfectness do  not pass to a subspace.

We briefly discuss the result $\dim_{\on{CA}}(X,d)  =\inf \{ \dim_{\on{A}}(X,\theta): \theta \in \sJ_p(X,d)\}$.
 Since $\sJ_p(X,d) \subset \sJ(X,d)$, the upper bound on $ \dim_{\on{CA}}(X,d)$ is obvious but the other inequality is non-trivial as it is possible that  $ \sJ_p(X,d) \neq \sJ(X,d)$ as we recall below. Trotsenko and  V\"ais\"al\"a characterize metric spaces for which  $ \sJ_p(X,d)= \sJ(X,d)$. To state their characterization, we recall the notion of \emph{weakly uniformly perfect} spaces. We say that a metric space is \emph{weakly uniformly perfect} if there exists $C>1$ such that if $B(x,r) \neq X$ for some $x \in X, r >0$, then either $B(x,r) = \{x\}$ or $B(x,r) \setminus B(x,r/C) \neq \emptyset$. The Trotsenko-V\"ais\"al\"a theorem states that a compact metric space $(X,d)$ satisfies  $ \sJ_p(X,d)= \sJ(X,d)$ if and only if $(X,d)$ is weakly uniformly perfect \cite[Theorems 4.11 and 6.20]{TV}.

\subsection{Outline of the work}
%The main result of this work is that the conformal Assouad dimension coincides with the critical exponent corresponding to the combinatorial modulus for any compact doubling metric space; that is,
%\be  \label{e:main}
%\dim_{\on{CA}}(X,d)= \on{CE}(X,d).
%\ee
% Our work builds upon earlier works of Carrasco  \cite{Car},   Keith  and Laakso \cite[\textsection 5.2]{KL}. 

 To show the estimate  $\dim_{\on{CA}}(X,d)\le \on{CE}(X,d)$, we construct a graph which is Gromov hyperbolic  called the hyperbolic filling (see \textsection \ref{s:construction}). A theorem of Bonk and Schramm implies that a quasi-isometric change of metric on the hyperbolic filling induces a power quasisymmetric change of metric on its boundary. 
Roughly speaking, a quasi-isometric change of metric is done using the optimal functions for the combinatorial modulus. This is done in \cite[Theorems 1.1 and 1.2]{Car} where the author introduces hypotheses on weight functions on the graph that defines a bi-Lipschitz change of metric in the hyperbolic filling. However the hypotheses   introduced in \cite[Theorem 1.1]{Car} implies that $(X,d)$ is uniformly perfect as pointed in \cite[Lemma 6.2]{Sha}. Since the metric spaces we consider are not necessarily   uniformly perfect, we need modify one of the hypothesis so that it is more suitable for bounding the conformal Assouad dimension (see hypothesis \hyperlink{H4}{$\on{(H4)}$} in Theorem \ref{t:H1-4}). The key new tool is a modification of a lemma of Vol'berg and Konyagin to construct a $p$-homogeneous measure  on $(X,\theta)$, where $\theta$ is power quasisymmetric to $d$ and $p>\on{CE}(X,d)$ (see Lemma \ref{l:redistribute} and Proposition \ref{p:assouad}).  This along with Theorem \ref{t:VK} implies the bound $\dim_{\on{CA}}(X,d)\le \on{CE}(X,d)$.
Another key distinction from \cite{Car} is that  the metric space is not necessarily uniformly perfect. Therefore by the Trotsenko-V\"ais\"al\"a theorem, this approach need not construct all possible metrics in $\sJ(X,d)$. Nevertheless, this approach  provides the sharp upper bound and also leads to  $\dim_{\on{CA}}(X,d)  =\inf \{ \dim_{\on{A}}(X,\theta): \theta \in \sJ_p(X,d)\}$.

 For the other bound $\on{CE}(X,d) \le \dim_{\on{CA}}(X,d)$, we use a $p$-homogeneous measure $\mu$ in $(X,\theta)$ and $\theta \in \sJ(X,d)$ for $p> \dim_{\on{CA}}(X,d)$ and define an function $\rho$ for the combinatorial modulus that is  similar to \cite[(3.7)]{Car}. However some modifications are needed because \cite{Car} uses the uniform perfectness in an essential way to control $\rho$. 
Some of the parameters and constants in \cite{Car} depend on the constant associated with the uniform perfectness property. Much of the work is about removing  such dependence on  uniform perfectness.

\section{Hyperbolic filling of a compact metric space}

\subsection{Gromov hyperbolic spaces and its boundary}
Let $(Z,d)$ be a metric space. 
	We recall some basic notions regarding Gromov hyperbolic spaces and refer the reader to \cite{BH, CDP,GH90,Gro87,Vai} for a detailed exposition. Given three points $x,y,w \in Z$, we define the \emph{Gromov product} of $x$ and $y$ with respect to the base point $w$ as
\[
(x|y)_w= \frac{1}{2} (d(x,w)+d(y,w)-d(x,y)).
\]
By the triangle inequality, Gromov product is always non-negative.
We say that a metric space $(Z,d)$ is $\delta$-\emph{hyperbolic}, if for any four points $x,y,z,w \in Z$, we have
\[
(x|z)_w \ge (x|y)_w \wedge (y|z)_w - \delta.
\]
We say that  $(Z,d)$ is hyperbolic (or $d$ is a hyperbolic metric), if $(Z,d)$ is hyperbolic for some $\delta \in [0,\infty)$.
If the above condition is satisfied for a fixed base point $w \in Z$, and arbitrary $x,y,z \in Z$, then $(Z,d)$ is $2 \delta$-hyperbolic \cite[Proposition 1.2]{CDP}.

We recall the definition of the boundary of a hyperbolic space. Let $(Z,d)$ be a hyperbolic space and let $w \in Z$.
A sequence of points $(x_i)_{i \in \bN} \subset Z$ is said to \emph{converge at infinity}, if
\[
\lim_{i,j \to \infty} (x_i|x_j)_w = \infty.
\] 
The above notion of convergence at infinity does not depend on the choice of the base point $w \in Z$, because by the triangle inequality $\abs{(x|y)_w- (x|y)_{w'}} \le d(w,w')$.

Two sequences $(x_i)_{i \in \bN}, (y_i)_{i \in \bN}$ that converge at infinity are said to be \emph{equivalent}, if 
\[
\lim_{i \to \infty} (x_i|y_i)_w = \infty.
\]
This defines an equivalence relation among all sequences that converge at infinity \cite[\textsection 1, Chapter 2]{CDP}.
As before, is easy to check that the notion of equivalent sequences does not depend on the choice of the base point $w$.
The \emph{boundary} $\partial Z$ of $(Z,d)$
is defined as the set of equivalence classes of sequences 
converging at infinity under the above equivalence relation.
If there are multiple hyperbolic metrics on the same set $Z$, to avoid confusion, we denote the boundary of $(Z,d)$ by $\partial(Z,d)$.
The notion of Gromov product can be defined on the boundary as follows: for all $a, b \in \partial Z$
\[
(a|b)_w = \sup \set{ \liminf_{i \to \infty} (x_i|y_i)_w : (x_i)_{i \in \bN} \in a, (y_i)_{i \in \bN} \in b}.
\]
By \cite[Remarque 8, Chapitre 7]{GH90}, if $(x_i)_{i \in \bN} \in a, (y_i)_{i \in \bN} \in b$, we have
\[
(a|b)_w -2 \delta \le \liminf_{i \to \infty} (x_i|y_i)_w \le (a|b)_w.
\]
%and similarly, for $a \in \partial Z, y \in Z$, we define
%\[
%(a|y)_w= \sup \set{ \liminf_{i \to \infty} (x_i|y)_w : \set{x_i} \in a}.
%\]
The boundary $\partial Z$ of the hyperbolic space $(Z,d)$ carries a family of metrics. A metric $\rho:\partial Z \times \partial Z \to [0,\infty)$ on $\partial Z$ is said to be a \emph{visual metric} with base point $w \in Z$ and \emph{visual parameter} $\alpha \in (1,\infty)$  if there exists $k_1,k_2 >0$ such that
\[
k_1 \alpha^{-(a|b)_w} \le \rho(a,b) \le k_2 \alpha^{-(a|b)_w} 
\]
 If a visual metric with base point $w$ and visual parameter $\alpha$ exists, then it can be chosen to be
\[
\rho_{\alpha,w}(a,b):= \inf \sum_{i=1}^{n-1}\alpha^{-(a_i|a_{i+1})_w},
\]
where the infimum is over all finite sequences $(a_i)_{i=1}^n \subset \partial Z, n \ge 2$ such that $a_1=a, a_n=b$. Any other visual metric with the same basepoint and visual parameter is bi-Lipschitz equivalent to $\rho_{\alpha,w}$.

Visual metrics exist on hyperbolic metric spaces as we recall now. % A metric space $(Z,d)$ is said to be \emph{proper} if all closed balls are compact.
For any  $\delta$-hyperbolic space $(Z,d)$, there exists $\alpha_0 >1$ ($\alpha_0$ depends only on $\delta$) such that if $\alpha \in (1,\alpha_0)$, then there exists a visual metric with parameter $\alpha$ \cite[Chapitre 7]{GH90}, \cite[Lemma 6.1]{BS}. %If a visual metric with visual parameter $\alpha$ exists, then it can be obtained (up to bi-Lipschitz equivlance)
It is well-known that quasi-isometry between almost geodesic hyperbolic spaces induces a quasisymmetry on their boundaries (the notion of almost geodesic space is given in Definition \ref{d:almostgeo}). Since this plays a central role in our construction of metric, we recall the relevant definitions and results below.

We say that a map (not necessarily continuous) $f:(X_1,d_1) \to (X_2,d_2)$ between two metric spaces is a \emph{quasi-isometry} if there exists constants $A,B>0$ such that
\[
A^{-1} d_1(x,y) - A \le d_2(f(x),f(y)) \le A d_1(x,y) +A
\]
for all $x,y \in X_1$, and
\[
\sup_{x_2 \in X_2} d(x_2,f(X_1)) = \sup_{x_2 \in X_2} \inf_{x_1 \in X_1} d(x_2,f(x_1)) \le B.
\]
	\begin{definition} \label{d:almostgeo}
	A metric space $(Z,d)$ is $k$-\emph{almost geodesic}, if
	for every $x, y \in Z$ and every $t \in [0, d(x,y)]$, there is some $z \in Z$ with
	$\abs{d(x,z) - t} \le k$ and $\abs{d(y,z) - (d(x,y) - t)} \le k$. We say that a metric space is \emph{almost geodesic} if it is $k$-almost geodesic for some $k\ge 0$. 
\end{definition}
Quasi-isometry between hyperbolic spaces induce quasisymmetries on their corresponding boundaries. We recall a result due to Bonk and Schramm below.
	\begin{proposition}[{\cite[Theorem 6.5 and Proposition 6.3]{BS}}] \label{p:BS}
	Let $(Z_1,d_1)$ and $(Z_2,d_2)$ be two almost geodesic, $\delta$-hyperbolic metric spaces. Let $f:(Z_1,d_1) \to (Z_2,d_2)$ be a  quasi-isometry. 
	\begin{enumerate}[(a)]
		\item If $(x_i)_{i \in \bN} \subset Z_1$ converges at infinity, then $(f(x_i))_{i \in \bN} \subset Y$ converges at infinity.
		If  $(x_i)_{i \in \bN}$ and $(y_i)_{i \in \bN}$ are equivalent sequences in $X$ converging at infinity, then
		$(f(x_i))_{i \in \bN}$ and $(f(y_i))_{i \in \bN}$ are also equivalent.
		\item  The map $\partial f: \partial Z_1 \to \partial Z_2$ given by $\partial f \left((x_i)_{i \in \bN}\right)= (f(x_i))_{i \in \bN}$ is well-defined, and is a bijection.
		\item Let $p_1 \in Z_1$ be a base point in $Z_1$, and let $f(p_1)$ be a corresponding base point in $Z_2$. Let $\rho_1, \rho_2$ denote  visual metrics (with not necessarily the same visual parameter) on $\partial Z_1, \partial Z_2$ with base points $p_1, f(p_1)$ respectively. Then the induced boundary map $\partial f: (\partial Z_1, \rho_1) \to (\partial Z_2, \rho_2)$ is a power quasisymmetry. 
	\end{enumerate}
\end{proposition}

\subsection{Geodesic hyperbolic spaces} \label{s:geod}
Let  $(Z,d)$ by a geodesic $\delta$-hyperbolic metric space. Recall that $(Z,d)$ is geodesic if for any $x,y \in X$, there exists a curve $\gamma:[0,d(x,y)] \to Z$ such that $\gamma(0)=x, \gamma(d(x,y))=y$ and $d(\gamma(s),\gamma(t))=\abs{s-t}$ for all $s,t \in [0,d(x,y)]$. Such a curve is called a geodesic between $x$ and $y$.
For $x,y \in Z$, we denote by $[x,y]$ a \emph{geodesic} between $x$ and $y$. For $x,y,z \in Z$, we denote by $[x,y,z]=[x,y] \cup [y,z] \cup [z,x]$ a \emph{geodesic triangle} in $Z$.  Recall that a \emph{tripod} is a metric tree with three edges arising from a common \emph{central vertex} such that each edge $a$ is isometric to the closed interval $[0,l(a)]$ for some $l(a) \ge 0$ called the length of the edge $a$.  A tripod is determined up to isometry by the length of the three edges. We allow for the degenerate case where the length of some of the edges could be zero.

 Given a geodesic triangle $\Delta=[x,y,z]$, there exists a map $f_\Delta: \Delta \to T_\Delta$ from $\Delta$ to a tripod $T_\Delta$ such that the restriction of $f_{\Delta}$ to each side of the triangle is an isometry  \cite[Proposition 2]{GH90}.  The \emph{inscribed triple} of a geodesic triangle $\Delta$ is defined to be the preimages of the `central vertex' of the tripod $T_\Delta$ under the map $f_{\Delta}$ described above.

 Unlike a tripod, a geodesic triangle $\Delta$ need not have a canonical center. However, it has a reasonable notion of approximate center. % We introduce the notion of a $K$-approximate center of a geodesic triangle.
  For $K \ge 0$, a point $c \in Z$ is a \emph{$K$-approximate center} of a geodesic triangle $[x,y,z]$ if $c$ is at a distance at most $K$ from each of the three sides, that is, $d(c,[x,y]) \vee d(c,[y,z]) \vee d(c,[z,x]) \le K$. The following proposition concerns a few properties of approximate center.
 
 \begin{proposition} \label{p:acent}
Let $(Z,d)$ be a geodesic, $\delta$-hyperbolic metric space. 
\begin{enumerate}[(a)]
	\item  Each point of the inscribed triple of a geodesic triangle is a $4 \delta$-approximate center.
	\item Any two $K$-approximate centers $c$ and $c'$ of a geodesic triangle $[x,y,z]$ satisfies $d(c,c') \le 8 K$.
	\item If $c$ is a $K$-approximate center of a geodesic triangle $[x,y,z]$ then
	\[
	\abs{d(x,c) - (y|z)_x} \le 4K.
	\]
	\item If $f:(Z_1,d_1) \to (Z_2,d_2)$ is a quasi-isometry between two geodesic $\delta$-hyperbolic metric spaces and if $c$ is a $K_1$-approximate center of $[x,y,z]$ then $f(c)$ is a $K_2$-approximate center of any geodesic triangle $[f(x),f(y),f(z)]$, where $K_2$ depends only on $K_1, \delta$ and the constants associated with the quasi-isometry $f$. In particular, 
	\[
	\abs{d_2(f(x),f(c)) - (f(y)|f(z))_{f(x)}} \le 4K_2.
	\]
\end{enumerate}
 \end{proposition}
\begin{proof}
	\begin{enumerate}[(a)]
		\item[(a)] By \cite[Proposition 21, Chapitre 2]{GH90} each point of the inscribed triple is a $4\delta$-approximate center.
		\item[(b,c)] Let $c$ denote a $K$-approximate center of $[x,y,z]$. Let $p_1,p_2,p_3$ be the points of the inscribed triple on $[x,y],[y,z],[z,x]$ respectively. Similarly, let $q_1,q_2,q_3$ be three points on $[x,y],[y,z],[z,x]$ respectively such that  $d(c,q_i) \le A$ for all $i=1,2,3$. This implies that $d(q_i,q_j) \le 2K$ for all $i,j$. By the argument in \cite[Proof of Lemme 20, Chapitre 2]{GH90} we have
		\be \label{e:ac1}
		d(p_i,q_i) \le 3K \quad \mbox{for all $i=1,2,3$.}
		\ee
		Since $d(x,p_1)=(y|z)_x$ and $d(p_1,q_1) \le 3K$, we obtain
		\[
		\abs{d(x,c) - (y|z)_x}= \abs{d(x,c)-d(x,p_1)} \le d(p_1,q_1)+d(c,q_1) \le 3K+K=4K.
		\]
		This concludes the proof of (c).
Similarly, $d(c,p_i) \le d(p_i,q_i) + d(c,q_i) \le 4K$ for all $i=1,2,3$. Therefore $d(c,c') \le d(c,p_1)+d(c',p_1) \le 8K$, and hence (b) holds.
\item[(d)]   This is an immediate consequence of the geodesic stability under quasi-isometries \cite[Th\'eor\`eme 11, Chapitre 5]{GH90} and (c).
	\end{enumerate}
\end{proof}
 
\subsection{Construction of hyperbolic filling} \label{s:construction}
In this section, we recall the construction of a hyperbolic filling of a compact   metric space.
Let $(X,d)$ be a compact metric space.
The construction below is due to A.~Bjorn, J.~Bjorn and Shanmugalingam \cite{BBS}. Earlier versions of this construction are due to Elek, Bourdon and Pajot \cite{Ele,BP}.

Let $\lambda,a \in (1,\infty)$ be two parameters which we call the horizontal and vertical parameter of the hyperbolic filling respectively.
We assume that the diameter is normalized so that $\diam(X,d) =\frac 1 2$. 
Let $X_n, n \in \bN_{\ge 0}$ be an increasing sequence of maximal $a^{-n}$-separated subsets of $X$. In other words, $X_n \subset X_m$ for all $n<m$,  any two distinct points in $X_n$ have mutual distance at least $a^{-n}$ and any set strictly larger than $X_n$ has two distinct points whose distance is strictly less than $a^{-n}$. The vertex set of the graph is $\sS= \cup_{n \ge 0} \sS_n$, where $\sS_n= \{(x,n): x \in X_n\}$. Two distinct vertices $(x,n),(y,m) \in \sS$ are joined by an edge if and only if either $n=m$ and $B(x, \lambda a^{-n}) \cap B(y, \lambda a^{-m}) \neq \emptyset$ or if $\abs{n-m}=1$ and $B(x,a^{-n}) \cap B(y,a^{-m}) \neq \emptyset$. 
 Let $D_1$ denote  the combinatorial (graph) distance on   $\sS$ defined by the above set of edges. That is $D_1((x,n),(y,m))$ is the minimal number $k$ such that $(x,n)=(x_0,n_0),(x_1,n_1),\ldots,(x_k,n_k)=(y,m)$, where $(x_i,n_i),(x_{i+1},n_{i+1}) \in \sS$ is joined by an edge for all $i=0,1,\ldots,k-1$. 
 It is evident that $(\sS,D_1)$ is $1$-almost geodesic metric space.

 We now construct a graph with fewer vertical edges. For each $(x,n) \in \sS$, we choose $(y,n-1) \in \sS_n$ such that $d(x,y) = \min_{z \in X_{n-1}} d(x,z)$. In this case, we say that $(y,n-1)$ is the \emph{parent} of $(x,n)$ or equivalently, $(x,n)$ is a \emph{child} of $(y,n-1)$. Such a choice of $y \in X_{n-1}$ is not unique but we fix this choice for the remainder of this work.  Since $X_{n-1}$ is a maximal $a^{-(n-1)}$ separated subset of $X$, $d(x,y) \le  a^{-(n-1)}$. Hence $x \in B(x,a^{-n}) \cap B(y,a^{-(n-1)}) \neq \emptyset$. In other words every parent and their child is connected by an edge in the graph associated with $(\sS,D_1)$. We define a new graph whose edges consists of all of the edges between parent and child and those between  $(x,n), (y,n) \in \sS$ where $B(x,\lambda a^{-n}) \cap B(y,\lam a^{-n}) \neq \emptyset, x \neq y$. The corresponding graph distance is denoted by $D_2$. The \emph{set of children of a vertex $v$} is denoted by 
 \be \label{e:defchild}
 C(v) := \{ w \in \sS : \mbox{$w$ is a child of $v$}\}.
 \ee
Note that $C(v) \subset \sS_{n+1}$ whenever $v\in \sS_n$.
 
 If $(x,n+1)$ and $(y,n+1)$ share an edge in $D_2$ and if $(x_0,n)$ and $(y_0,n)$ are their respective parents, then
 $d(x_0,y_0) \le d(x,y) +d(x,x_0)+d(y,y_0) < 2 a^{-n}+ 2 \lambda a^{-n-1}$. 
Under the assumption $\lam \ge 2+ 2 \lam a^{-1}$, we have $D_2((x,n+1),(y,n+1)) \le 1$, then $D_2((x_0,n),(y_0,n)) \le 1$
 whenever $(x_0,n),(y_0,n)$ are the parents of $(x,n+1),(y,n+1)$ respectively. We say that $(x,n)$ is a \emph{descendant} of $(y,k)$ if $n>k$, and there exists $(z_j,n_j) \in \sS$ for $j=0,\ldots,n-k$ such that $(z_0,n_0)=(y,k), (z_{n-k},n_{n-k})=(x,n)$, where $(z_{i+1},n_{i+1})$ is a child of $(z_i,n_i)$ for all $i=0,\ldots,n-k-1$. For any $n > k \ge 0$ and $v \in \sS_k$, the \emph{set of descendants of $v$ in generation $n$} is denoted by 
 \be \label{e:defdes}
 \sD_n(v) = \{w \in \sS_n: \mbox{$w$ is a descendant of $v$}\}.
 \ee
The following lemma 
  is an analogue of \cite[Lemma 2.2]{Car}.
\begin{lem} \label{l:filling}
	Let $\lam,a >1$ be horizontal and vertical parameters of the hyperbolic filling respectively. 
	\begin{enumerate}[(a)]
		\item If $(z,n+1)$ is a child of $(x,n)$, then $d(x,z) < a^{-n}$. If $(y,k)$ is a descendant of $(x,n)$ (for some $k>n$), then 
		\[
		d(x,y) < \frac{a}{a-1} a^{-n}.
		\]
		
		\item If $\lam \ge 2+2 \lam a^{-1}$ and
		$D_2((x,n+1),(y,n+1)) \le 1$, then $D_2((x_0,n),(y_0,n)) \le 1$, where $(x_0,n),(y_0,n)$ are the parents of $(x,n+1), (y,n+1)$ respectively.  Similarly, if  $\lam \ge 2+4 \lam a^{-1}$ and
		$D_2((x,n+1),(y,n+1)) \le 2$, then $D_2((x_0,n),(y_0,n)) \le 1$, where $(x_0,n),(y_0,n)$ are the parents of $(x,n+1), (y,n+1)$.
		\item Let $\lam \ge 6$. If $(x,n+1),(y,n+1) \in \sS_{n+1}$ such that $d(x,y) \le 4 a^{-n}$. If $(x_0,n), (y_0,n) \in \sS_n$ are the parents of  $(x,n+1),(y,n+1)$ respectively, then $D_2((x_0,n),(y_0,n))\le 1$.
		\item If $\lam>1 +a^{-1}$, we have $D_1 \le D_2 \le 2D_1$.
		\item Let $\lam > 1 +a^{-1}$. Let $w \in \sS_{n+1}$ and $u,v \in \sS_n$ be such that $D_1(u,w)=1$ and $D_2(v,w)=1$. Then $D_2(u,v) \le 1$.
	\end{enumerate}
\end{lem}
\begin{proof}
	\begin{enumerate}[(a)]
		\item Since $X_n$ is maximal $a^{-n}$-subset of $X$, every point $z \in X$ satisfies $d(z,X_n) < a^{-n}$. This shows the first claim. 
		If $(y,k)$ is a descendant of $(x,n)$ by the first claim and triangle inequality $d(x,y) \le \sum_{i=n}^k a^{-i} < \frac{a}{a-1} a^{-n}$.
		\item  Since
		$d(x_0,y_0) \le d(x,y) +d(x,x_0)+d(y,y_0) < 2 a^{-n}+ 2 \lambda a^{-n-1}$,  we have $ \{x_0,y_0\} \subset B(x_0,\lam a^{-n}) \cap B(y_0, \lam a^{-n}) \neq \emptyset$ for any $\lambda$ such that $\lambda \ge 2 + 2 \lambda a^{-1}$. Hence  $D_2((x_0,n),(y_0,n)) \le 1$. The other claim follows from a similar argument.
		\item Since $d(x_0,y_0) \le d(x,y) +d(x,x_0)+d(y,y_0) < 4a^{-n} + a^{-n} + a^{-n}=6a^{-n}$, we have $ \{x_0,y_0\} \subset B(x_0,\lam a^{-n}) \cap B(y_0, \lam a^{-n}) \neq \emptyset$ whenever $\lam \ge 6$. Therefore $D_2((x_0,n), (y_0,n)) \le 1$.
		\item [(d),(e)]
		 Since every edge in the graph corresponding to $(\sS,D_2)$ is contained in the graph corresponding to $(\sS,D_2)$, we have $D_2 \ge D_1$.

		On the other hand, if there is an edge in $(\sS,D_1)$\footnote{Here we abuse notation and use $(\sS,D_i)$ to denote the graph, for $i=1,2$.} which is not present in $(\sS,D_2)$, then it must be  between  some $(x,n), (y,n+1)$ such that $n \in \bN_{\ge 0}$ and that the parent of $(y,n+1)$ is $(z,n)$ where $z \neq x$. In this case $d(x,y) \le   a^{-n-1}+a^{-n}$ (since $B(x,a^{-n}) \cap B(y,a^{-n-1}) \neq \emptyset$). Therefore if $\lambda > 1+ a^{-1}$, there would be an edge between $(x,n)$ and $(z,n)$ in both graphs (since $y \in B(x, \lam a^{-n}) \cap B(z, \lam a^{-n})$). This implies that
		\be 
		D_1 \le D_2 \le 2 D_1, \quad \mbox{whenever $\lam\ge 1+a^{-1}$.}
		\ee

	\end{enumerate}
\end{proof}

 We recall the relevant properties the metric spaces $(\sS,D_1)$ and $(\sS,D_2)$. By the choice of the $\diam (X,d)$, there is an unique point $x_0 \in X_0$. We choose $v_0:=(x_0,0)$ as the base point of the metric spaces $(\sS,D_1)$ and $(\sS,D_2)$. We denote the Gromov product with respect to the basepoint $v_0$ in $(\sS,D_1)$ and $(\sS,D_2)$ by $(\cdot|\cdot)_1, (\cdot | \cdot)_2$ respectively. The key point in the following result is that the hyperbolicity constant $\delta$ depends only on $a$ and $\lam$ unlike the analogous result in \cite{BP,Car} where $\delta$ also depends on the constant associated with the uniform perfectness property (see \cite[Remark after Proposition 2.1]{Car}).
 \begin{prop} (Cf. \cite[Lemma 3.3 and Theorem 3.4]{BBS}) \label{p:bbs}
 	Let $(X,d)$ be a compact metric space and let $a,\lambda$ denote the vertical and horizontal parameters respectively of the hyperbolic filling. Then the hyperbolic filling $(\sS,D_1)$ satisfies the following properties
 	\begin{enumerate}[(a)]
 		\item For any $v=(z,n),w=(y,m) \in \sS$, we have
 		\[
 	\frac{a-1}{4 \lam a} \left(d(z,y)+a^{-n}+a^{-m}\right)	\le a^{-(v|w)_1} \le \frac{a^{5/2}}{\lam -1} \left(d(z,y)+a^{-n}+a^{-m}\right)
 		\]
 		 In particular, if $a,\lam \in [2, \infty)$ and $a \ge \lam$, then
 		 \[
 		 \abs{(v|w)_1+ \frac{\log(d(x,y)+a^{-m}+a^{-n})}{\log a}} \le 4.
 		 \]
 		\item $(\sS,D_1)$ is $\delta$-hyperbolic, where $\delta=2 \frac{\log \left(\frac{8 \lam a^{7/2}}{(a-1)(\lam -1)}\right)}{\log a}$. In particular, if $a,\lam \in [2, \infty)$ and $a \ge \lam$ implies that $\delta$ can be chosen to be $15$.
 	\end{enumerate}
 \end{prop}
\begin{proof}
	\begin{enumerate}[(a)]
		\item  The first estimate follows from \cite[Proof of Lemma 3.3]{BBS}. The second conclusion is a consequence of the estimate
		\[
	 \max \left( \frac{\log(a^{5/2}/(\lam -1))} {\log a}, \frac{\log \left(\frac{4 \lam a}{a-1}\right)}{\log a} \right) \le 4 \quad \mbox{whenever $a \ge \lam \ge 2$.}
		\]
		\item The $\delta$-hyperbolicity follows from the proof of \cite[Theorem 3.4]{BBS} along with \cite[Proposition 1.2]{CDP}.
For the second conclusion, observe that
\[
2 \frac{\log \left(\frac{8 \lam a^{7/2}}{(a-1)(\lam -1)}\right)}{\log a} \le 2 \frac{\log (32a^{5/2})}{\log a}\le 15, \quad \mbox{whenever $a \ge \lam \ge 2$.}
\]
	\end{enumerate}
\end{proof}

By Proposition \ref{p:bbs}(a),
 a sequence of vertices $ ((x_i,n_i))_{i \in \bN} \in \sS$ converges at infinity if and only if $\lim n_i =\infty$ and $(x_i)_{i \in \bN}$ is a convergent sequence in $(X,d)$. Two sequences $((x_i,n_i))_{i \in \bN}$ and $((y_i,m_i))_{i \in \bN}$ that converge at infinity are equivalent if and only if $\lim_{i \to \infty} x_i= \lim_{i \to \infty} y_i$ and $\lim_{i \to \infty} n_i=\lim_{i \to \infty} m_i =\infty$. 
Define the \emph{limit at infinity} function $\mathfrak{l}_\infty: \partial (\sS,D_1) \to X$ that maps an equivalence class of sequence converging at infinity to its \emph{limit} as
\begin{equation} \label{e:deflimit}
\mathfrak{l}_\infty( ((x_i,n_i))_{i \in \bN})=\lim_{i \to \infty} x_i.
\end{equation}
 Note that $\mathfrak{l}_\infty$ is well defined and is a bijection. By the first estimate in Proposition \ref{p:bbs}(a),  $d$ is a visual metric on $\partial (\sS,D_1)$ with visual parameter $a$ and base point $v_0$ in the following sense:
\begin{equation} \label{e:lvisual}
\frac{a-1}{4 \lam a} d(x,y)	\le a^{-(\mathfrak{l}_\infty^{-1}(x)|\mathfrak{l}_\infty^{-1}(y))_1} \le \frac{a^{5/2}}{\lam -1} d(x,y) \quad \mbox{for all $x,y \in X$.}
\end{equation}

We would like to use Proposition \ref{p:acent} to estimate Gromov product in the hyperbolic filling. Since $(\sS,D_1)$ is not a geodesic space, we embed it into a geodesic space by replacing each edge with an isometric copy of the unit interval to obtain a metric space $(\wt{\sS},D_1)$ where we view $\sS \subset \wt \sS$ and $D_1$ on $\wt \sS$ is an extension of $D_1$ on $\sS$. 
 For any $x,y \in X$, let $\wt n$ is the largest integer that satisfies $\{x,y\} \subset B(\wt z,2 a^{-\wt n})$ for some $(\wt z,\wt n) \in \sS$ and define
 \be \label{e:defcxy}
 c(x,y)= \{(\wt z,\wt n) \in \sS: \{x,y\} \subset B(\wt z,2 a^{-\wt n}) \}.
 \ee
 We think of $c(x,y)$ as the set of approximate centers of the   triangle 
$[v_0, \mathfrak{l}_\infty^{-1}(x), \mathfrak{l}_\infty^{-1}(y)]$, where $\mathfrak{l}_\infty^{-1}(x),\mathfrak{l}_\infty^{-1}(y) \in \partial (\mathcal{S},D_1)$ and $v_0$ is the unique element of  $\mathcal{S}_0$.
The following lemma makes this precise by identifying $c(x,y)$ as approximate centers of certain geodesic triangles in $\wt \sS$.

\begin{lem} \label{l:centhyp}
	Let the parameters of the hyperbolic filling satisfy
	\[
	a \ge \lam \ge 2, \quad \mbox{and} \quad \lam \ge 1+\frac{a}{a-1}.
	\]		
 Let $(z,m),(w,n) \in \sS$ such that $x \in B(z,a^{-m}), y\in B(w,a^{-n})$ and such that $a^{-m}+a^{-n} < a^{-2} d(x,y)$. Then any $(w,k) \in c(x,y)$ is a $K$-approximate center for any geodesic triangle $[v_0,(z,m),(w,n)]$ in $(\wt \sS,D_1)$, where $K=80$.
\end{lem}
\begin{proof}
	Since every point in $\wt \sS$ is at most distance $\frac{1}{2}$ away from a point in $\sS$, by replacing points in $\wt \sS$ with the corresponding closest points in $\sS$, we obtain that $(\wt \sS,D_1)$ is $(\delta+3)$-hyperbolic whenever $(\sS,D_1)$ is $\delta$-hyperbolic.
	
	Let $(\wt z,k) \in c(x,y)$. 
	Since $\{x,y\} \subset B(\wt z,2a^{-k})$, we obtain $d(x,y) \le 4 a^{-k}$. Choose $w \in X$ such that  $(w,k+1) \in \sS_{k+1}$ and $d(x,w)<a^{-(k+1)}$. By the maximality of $k$,  we have that $y \notin B(w,2a^{-(k+1)})$ and hence $d(x,y) \ge d(w,y) - d(w,x)>a^{-(k+1)}$. In particular,
	\be \label{e:boundm}
	a^{-(k+1)} < d(x,y) \le 4 a^{-k}.
	\ee
	If $a \ge 2$, we have
	$a^{-(k+1)} < d(x,y) \le  4 a^{-n}\le a^{-k+2}$, which implies 
	\be \label{e:hyp1}
	-1 \le k + \frac{\log d(x,y)}{\log a} \le 2, \quad \mbox{whenever $a \ge 2$.}
	\ee
	Since $a^{-m} <a^{-2}d(x,y)$, we have $m+\frac{\log d(x,y)}{\log a}>2$ which along with \eqref{e:hyp1} implies that $m \ge k$. Choose $(w,k)$ such that $(z,m)$ is a descendant of  $(w,k)$. By Lemma \ref{l:filling}(a),
	we have
	$d(x,w) \le d(x,z)+d(w,z) < a^{-m}+ \frac{a}{a-1} a^{-k} \le \left(1+\frac{a}{a-1}\right) a^{-k}$. Therefore if $\lam \ge 1+\frac{a}{a-1}$, we have $x \in B(w,\lam a^{-k}) \cap B(\wt z,\lam a^{-k}) \neq \emptyset$. Hence  $(w,k)$ and $(\wt z,k)$ are  either equal or horizontal neighbors.
	
	Note that 
	$\abs{d(x,y)-d(z,w)} \le a^{-m}+a^{-n} <a^{-2}d(x,y)$, which implies
		\[
(1-a^{-2})d(x,y)	\le 	d(z,w)+a^{-m}+a^{-n} \le (1+2a^{-2}) d(x,y). 
		\]
	Therefore if $a \ge \lam \ge 2$, 
	we have
	\[
	\abs{ \frac{\log(d(x,y))}{\log a}- \frac{\log(d(z,w)+a^{-m}+a^{-n})}{\log a}} \le 1.
	\]
	Combining with \eqref{e:hyp1} and Proposition \ref{p:bbs}(a), we obtain that
	\[
	\abs{((z,m)|(w,n))_1-k} \le 7.
	\]
	This along with Proposition \ref{p:acent}(a) and the fact that $(\wt z,k)$ is a neighbor of $(w,k)$, we obtain that 
	$(\wt z,k)$ is $K$-approximate center of the geodesic triangle $[v_0,(z,m),(w,n)]$, where $K=1+7+4(15+3)=80$.
\end{proof}
\begin{remark}
	The assumption 
	\be  \label{e:assume}
	a \ge \lam \ge 6
	\ee
	implies the estimates assumed on $a,\lam$ in Lemma \ref{l:filling}, Proposition \ref{p:bbs}, and Lemma \ref{l:centhyp} hold. For this reason we  assume \eqref{e:assume} for much of this work. The analogous estimate \cite[(2.8)]{Car} is more complicated  because it involves the constant in the definition of uniform perfectness.
\end{remark}
\section{Construction of metric  and homogeneous measure}

In this section, we construct metric in the conformal gauge and a homogeneous measure using a weight function on the hyperbolic filling $\sS$ as constructed in \textsection \ref{s:construction}. A \emph{weight} on a filling $\sS$ is   a function $\rho: \sS \to (0,\infty)$. We recall the definition of homogeneous measure and its relevance to Assouad dimension in \textsection \ref{s:VK}. We then introduce and recall some hypothesis on a weight function on the hyperbolic filling that provides upper bound on $\dim_{\on{CA}}(X,d)$ in \textsection \ref{s:hyp}. 

\subsection{Vol'berg-Konyagin theorem} \label{s:VK}
Our approach to obtain upper and lower bounds on the conformal Assouad dimension 
($\dim_{\on{CA}}(X,d) \le \on{CE}(X,d)$ and  $\dim_{\on{CA}}(X,d) \ge \on{CE}(X,d)$) relies on a theorem of Vol'berg and Konyagin that we recall below in Theorem \ref{t:VK}.
This result clarifies the relationship between Assouad dimension and  doubling measures. 
A non-zero Borel measure $\mu$ on a metric space $(X,d)$ is said to be \emph{doubling} if there exists $C_D>1$ such that 
\[
\mu(B(x,2r)) \le C_D \mu(B(x,r)) \quad \mbox{for all $x \in X, r>0$,}
\]
where $B(x,r)= \set{y \in X: d(x,y) <r }$ denotes the open ball of radius $r$ centered at $x$. A non-zero Borel measure is said to be \emph{$q$-homogeneous measure} if there exists $C>1$ such that 
\[
\mu(B(x,R)) \le C \left(\frac R r\right)^q \mu(B(x,r)), \quad \mbox{for all $x \in X, 0 < r  \le R$.}
\] 
It is evident that a measure is doubling if and only if it is $q$-homogeneous for some $q \in (0,\infty)$. The fundamental relationship between Asssouad dimension and doubling measures is given by the following theorem of Vol'berg and Konyagin \cite{VK}.
\begin{theorem} \label{t:VK} \cite[Theorem 1]{VK}  
	The Assouad dimension of compact metric space $(X,d)$ is given by
	\[
	\dim_{\on{A}}(X,d) = \inf  \{q >0 : \mbox{there exists a $q$-homogeneous measure on $(X,d)$}\}.
	\]
\end{theorem}

\subsection{Weights on the filling and Carrasco-type hypotheses} \label{s:hyp}

Let $\pi_1:\sS \to X, \pi_2:\sS \to \bN_{\ge 0}$ denote the projection maps such that $v=(\pi_1(v),\pi_2(v))$ for any $v \in \sS$. We say that an edge between two vertices $v$ and $w$ is \emph{horizontal} if $\pi_2(v)=\pi_2(w)$. For a vertex $v \in \sS$, by $B_v$ we denote the metric ball $B(\pi_1(v),a^{-\pi_2(v)})$.
Given a vertex $v \in \sS$, we define  the \emph{genealogy} $g(v)$ as a sequence of vertices $(v_0,v_1,\ldots,v_k)$ where $v_k=v$ and $v_i$ is the parent of $v_{i+1}$ for all $i=0,\ldots,k-1$ and $v_0$ is the unique vertex in $\sS_0$. 
Given a function $\rho:\sS\to (0,\infty)$, we define $\pi:\sS \to (0,\infty)$ as 
\be \label{e:defpi}
\pi(v)=\prod_{w \in g(v)} \rho(w).
\ee
A \emph{path} $\gamma$ in $(\sS,D_2)$ is a sequence of vertices $\gamma=(w_1,\ldots,w_n)$ where there is an edge between $w_i$ and $w_{i+1}$ (that is, $D_2(w_i,w_{i+1})=1$) for each $i=1,\ldots,n-1$. The $\rho$-length of a path $\gamma$ is defined by
\be  \label{e:defLrho}
L_\rho(\gam)=\sum_{v \in \gam} \pi(v).
\ee
The following two families of paths will play an important role in this work. A path is said to be \emph{horizontal} if it only consists of horizontal edges. Given $x,y \in X$ and $n \in \bN_{\ge 0}$, we define 
\be \label{e:defGamman}
\Gam_n(x,y)= \Biggl\{\gam=(v_1,\ldots,v_k) \Biggm|	\begin{minipage}{235 pt} $\gam$ is a path in $(\sS,D_2)$, $\pi_2(v_1)=\pi_2(v_k)=n$, and
	$x \in B_{v_1}, y \in B_{v_k}, k \in \bN$ \end{minipage}\Biggr\}.
\ee
For a vertex $v \in \sS_k$, we define 
\begin{equation} \label{e:defGamk}
	\Gam_k(v)=\inf \Biggl\{ \gam=(v_1,v_2,\ldots,v_n) \Biggm|
	\begin{minipage}{260 pt}
		$\gam$ is a horizontal path with $\pi_2(v_i)=k+1$ for all $i$, $n \in \bN$, $\pi_1(v_1) \in B_v, \pi_1(v_n) \notin 2\cdot B_v$	\end{minipage}
	\Biggr\}.
\end{equation}
For $x,y \in X$, we define
\be \label{e:defpic}
\pi(c(x,y))= \max\{\pi(w): w \in c(x,y)\},
\ee
where $c(x,y)$ is as defined in Lemma \ref{l:centhyp}.
We recall the Carrasco-type conditions imposed on the weight function $\rho: \sS \to (0,1)$. 

\begin{assumption} \label{a:h1-h4}
	A weight function $\rho: \sS \to (0,1)$ may satisfy some of the following hypotheses:
	\begin{enumerate}
		\item[ \hypertarget{H1}{$\on{(H1)}$}]  There exist $0<\eta_-\le \eta_+<1$ so that $\eta_- \le \rho(v) \le \eta_+$ for all $v \in \sS$.
		\item[ \hypertarget{H2}{$\on{(H2)}$}] There exists a constant $K_0\ge 1$ such that for all $v,w \in \sS$ that share a horizontal edge, we have
		\[
		\pi(v) \le K_0 \pi(w),
		\]
		where $\pi$ is as defined in \eqref{e:defpi}.
		\item[ \hypertarget{H3}{$\on{(H3)}$}]  
		There exists   a constant $K_1\ge 1$ such that for any pair of points $x,y\in X$, there exists $n_0 \ge 1$ such that if $n \ge n_0$ and $\gamma$ is a path in $\Gamma_n(x,y)$, then 
		\[
		L_\rho(\gamma) \ge K_1^{-1} \pi(c(x,y)),
		\]
		where $\Gamma_n(x,y),L_\rho, \pi(c(x,y))$ are as defined in \eqref{e:defGamman}, \eqref{e:defLrho}, and \eqref{e:defpic} respectively. 
		\item[ \hypertarget{H4}{$\on{(H4)}$}]   There exists $p>0$  such that for all $v \in \sS_m$ and $n >m$, we have
		\[
		\sum_{w \in \sD_n(v)} \pi(w)^{p} \le  \pi(v)^{p},
		\]
		where $\sD_n(v)$ denotes the descendants of $v$ in $\sS_n$. Clearly, it suffices to impose the above condition only for $n=m+1$.
		\item[ \hypertarget{H4-old}{$\on{(\wt{H4})}$}]   There exists $C>1, p>0$  such that for all $v \in \sS_m$, and $n >m$, we have
		\[
		C^{-1}  \pi(v)^{p} \le \sum_{w \in \sD_n(v)} \pi(w)^{p} \le C  \pi(v)^{p},
		\]
		where $\sD_n(v)$ denotes the descendants of $v$ in $\sS_n$.
	\end{enumerate}
\end{assumption}
The main results of this section are
Theorems \ref{t:H1-4} and \ref{t:simp} which provide an upper bound on $\dim_{\on{CA}}(X,d)$ under certain assumptions on the weight function on a filling.
Theorem \ref{t:H1-4} is an analogue of \cite[Theorem 1.1]{Car} but the hypothesis  \hyperlink{H4}{$\on{(H4)}$} is   different from that of \cite[Theorem 1.1]{Car} where the hypothesis \hyperlink{H4-old}{$\on{(\wt{H4})}$} was used instead. 
As explained in the introduction, \cite[Lemma 6.2]{Sha} implies that  \hyperlink{H1}{$\on{(H1)}$} along with  \hyperlink{H4-old}{$\on{(\wt{H4})}$}  can hold only on  a uniformly perfect metric space. Since we consider metric spaces that are not necessarily uniformly perfect, we need to modify \hyperlink{H4-old}{$\on{(\wt{H4})}$} to \hyperlink{H4}{$\on{(H4)}$}  as above. This hypothesis plays a key role in the upper bound on Assouad dimension. Another distinction from \cite{Car} is that the weights can be used to construct essentially all metrics in $\sJ(X,d)$. Since our construction of metric relies on Proposition \ref{p:BS}(c), we can only construct metrics in $\sJ_p(X,d)$ and hence we cannot obtain such a result. The following theorem is a counterpart of  \cite[Theorem 1.1]{Car}.
\begin{theorem} \label{t:H1-4}  
	Let $(X,d)$ be a compact doubling metric space and let $\sS$ be a hyperbolic filling with vertical and horizontal parameters $a, \lam$ such that $a \ge \lam \ge 6$. Let $\rho:\sS \to (0,1)$ be a weight function satisfying the hypothesis  \hyperlink{H1}{$\on{(H1)}$},  \hyperlink{H2}{$\on{(H2)}$},  \hyperlink{H3}{$\on{(H3)}$}, and  \hyperlink{H4}{$\on{(H4)}$}. Then there exists $\Theta_\rho \in \sJ_p(X,d)$ such that  $\dim_{\on{A}}(X,\Theta_\rho) \le p$.
\end{theorem}
 
The following is an analogue of \cite[Theorem 1.2]{Car}.
The conclusion of Theorem \ref{t:simp} is the same as that of Theorem \ref{t:H1-4} but the hypotheses  \hyperlink{H1}{$\on{(H1)}$},  \hyperlink{H2}{$\on{(H2)}$},  \hyperlink{H3}{$\on{(H3)}$}, and  \hyperlink{H4}{$\on{(H4)}$} are replaced by simpler hypotheses \hyperlink{S1}{$\on{(S1)}$} and \hyperlink{S2}{$\on{(S2)}$}.
 The hypotheses \hyperlink{S1}{$\on{(S1)}$} and \hyperlink{S2}{$\on{(S2)}$} below are identical to \cite[Theorem 1.2]{Car} but in the conclusion we bound $\dim_{\on{CA}}$ instead of $\dim_{\on{ARC}}$.
\begin{theorem} \label{t:simp}
	Let $(X,d)$ be a compact doubling metric space and let $\sS$ denote a hyperbolic filling with vertical and horizontal parameters $a,\lam$ respectively such that $a \ge \lam \ge 6$.
	Let $p>0$. There exists $\eta_0 \in (0,1)$ which depends only on $p, \lam$ and the doubling constant of $(X,d)$ (but not on the vertical parameter $a$) such that if there exists a function $\sigma: \sS \to [0,\infty)$   that satisfies:
	\item [\hypertarget{S1}{$\on{(S1)}$}]
	for all $v \in \sS_k$ and $k \ge 0$, if $\gam \in \Gam_k(v)$, then 
	\[
	\sum_{w \in \gam} \sigma(v) \ge 1,
	\]
	where $\Gam_k$ is as defined in \eqref{e:defGamk}, and
	\item [\hypertarget{S2}{$\on{(S2)}$}] 
	for all $k \ge 0$ and all $v \in \sS_k$, we have 
	\[
	\sum_{w \in C(v)} \sigma(v)^p \le \eta_0,
	\]
 then there exists $\Theta \in \sJ(X,d)$ such that $\dim_{\on{A}}(X,\Theta) \le p$ and $\Theta \in \sJ_p(X,d)$. In particular,  $\dim_{\on{A}}(X,d) \le \inf \{\dim_{\on{A}}(X,\theta): \theta \in \sJ_p(X,d)\} \le p$.
\end{theorem}

\subsection{Construction of metric using weights on the filling}
In this subsection, given a weight function which satisfies the hypotheses \hyperlink{H1}{$\on{(H1)}$}, \hyperlink{H2}{$\on{(H2)}$}, \hyperlink{H3}{$\on{(H3)}$}, we construct a metric  $\Theta_\rho \in \sJ_p(X,d)$ (Corollary \ref{c:visual}). The main idea is to use weight $\rho$ to induce a quasi-isometric change of metric on $\sS$ (Lemma \ref{l:Drho}(a)) which in turn induces a power quasisymmetric change of metric on its boundary by Proposition \ref{p:BS}(c). Since the boundary of $\sS$ can be identified with $(X,d)$ by Proposition \ref{p:bbs}(c), we therefore obtain a metric $\Theta_\rho \in \sJ_p(X,d)$. We remark that the hypothesis \hyperlink{H4}{$\on{(H4)}$} will not play any role in this subsection but will play a central role in the next one.

The weight function $\rho$ which satisfies the hypotheses \hyperlink{H1}{$\on{(H1)}$}-\hyperlink{H3}{$\on{(H3)}$} induces a metric $D_\rho$ on $\sS$. We set the length of an edge $e=(v,w)$ as $\ell_\rho(e)$, where
\[
\ell_\rho(e)= \begin{cases}
	2 \max \{-\log \eta_-,(-\log \eta_+)^{-1}, \log K_0\} & \mbox{if $e$ is horizontal},\\
	\log \frac{1}{\rho(v)} &\mbox{if $w$ is the parent of $v$,}
\end{cases}
\]
where $\eta_-,\eta_+,K_0$ are as defined in \hyperlink{H1}{$\on{(H1)}$} and \hyperlink{H2}{$\on{(H2)}$}.
This defines a metric
\be \label{e:defDrho}
D_\rho(v,w)= \inf_{\gamma} \sum_{e \in \gam} \ell_\rho(e),
\ee
where $\gamma$ varies over all paths in the graph $(\sS,D_2)$ from $v$ to $w$ and $e$ varies over all edges in $\gam$. By replacing each edge $e$ with an isometric copy of the interval $[0,\ell_\rho(e)]$, we define a geodesic metric space $(\wt \sS, D_\rho)$ such that $\sS \subset \wt \sS$ and the restriction of $D_\rho$ on $\wt \sS$ coincides with that of $\sS$.
\begin{lem}   \label{l:Drho}
	Let $(X,d)$ and $\sS$ be as in the statement of Theorem \ref{t:H1-4}.
	Let $\rho:\sS \to (0,\infty)$ which satisfies hypotheses  \hyperlink{H1}{$\on{(H1)}$} and \hyperlink{H1}{$\on{(H2)}$}.  Let $D_\rho$ denote the metric defined in \eqref{e:defDrho}.
	\begin{enumerate}[(a)]
		\item $(\sS,D_\rho)$ is approximately geodesic. The identity map $\on{Id}:(\sS,D_1) \to (\sS,D_\rho)$ is a quasi-isometry.

		\item Any path  of the form $((x_i,n_i))_{1\le i \le k}$ such that $(x_i,n_i)$ is the parent of $(x_{i+1},n_{i+1})$ for all $i=1,\ldots,k-1$ defines a shortest path in the $D_\rho$ metric and hence
		\[
		D_\rho((x_1,n_1),(x_k,n_k))= \abs{\log \frac{1}{\pi((x_1,n_1))} -\log \frac{1}{\pi((x_k,n_k))} }.
		\] 
		\item A sequence of vertices converges at infinity in $(\sS,D_1)$ if and only if it converges at infinity   $(\sS,D_\rho)$. Two sequences that converge at infinity are equivalent in $(\sS,D_1)$ if and only if they are equivalent in $(\sS,D_\rho)$. In particular, the identity map is a well-defined bijection between $\partial (\sS,D_1)$ and $\partial (\sS,D_\rho)$. 
		Therefore the limit at infinity map $\mathfrak{l}_\infty: \partial (\sS,D_1) \to X$ defined in \eqref{e:deflimit} is also well-defined as  $\mathfrak{l}_\infty:  \partial (\sS,D_\rho) \to X$. Furthermore, there exists
		 $C>1$ such that the Gromov product satisfies 
		\[
	 C^{-1} C \pi(c(x,y)) \le	e^{-(\mathfrak{l}_\infty^{-1}(x)|\mathfrak{l}_\infty^{-1}(y))_\rho} \le C \pi(c(x,y)) \quad \mbox{for all $x,y \in X$,}
		\]
		where $(\cdot | \cdot)_\rho$ is the Gromov product on $(\sS,D_\rho)$ with base point $v_0 \in \sS_0$.
		\item There exists $C>0$ such that the following holds: for any pair of distinct points $x,y \in X$, there exists $n_0$ such that whenever $n \ge n_0$ and $u,v \in \sS_n$ such that $x \in B_{u}, y\in B_v$, there exists a path $\gam=(w_i)_{i=0,\ldots,k}$ in the graph $(\sS,D_2)$ such that
		$L_\rho(\gam) \le C \pi(c(x,y))$.
	\end{enumerate}
\end{lem}
\begin{proof}
	\begin{enumerate}[(a)]
		\item  It is easy to check that $(\sS,D_\rho)$ is $2 \max \{-\log \eta_-,(-\log \eta_+)^{-1}, \log K_0\}$-approximately geodesic, since the horizontal edges are the longest edges. The fact that the identity map is a quasi-isometry is because there exist constants $C_1,C_2$ such that $C_1 \le \ell_\rho(e) \le C_2$ for all edges $e$. This along with Lemma \ref{l:filling}(d) implies that $C_1 D_1 \le C_1 D_2 	\le D_\rho \le C_2 D_2 \le 2 C_2 D_1$.
		\item This follows from the same argument as \cite[Proof of Lemma 2.3]{Car} which uses Lemma \ref{l:filling}(b).
		\item The first three claims follow from Proposition \ref{p:BS}. Let $(v_n)_{n \in \bN}, (w_n)_{n \in \bN}$ be two sequences
		of vertices such that $x \in B_{v_n}, y \in B_{w_n}, v_n \in \sS_n, w_n \in \sS_n$ for all $n \in \bN$. 
		By Lemma \ref{l:centhyp}, every vertex in $c(x,y)$ is a 80-approximate center for the geodesic triangle $[v_0,v_n,w_n]$
		(in $(\wt  \sS,D_1)$) for all large enough $n$. 
		By (a) and Proposition \ref{p:acent}(d), every vertex in $c(x,y)$ is a $K'$-approximate center of the geodesic triangle $[v_0,v_n,w_n]$ (in $(\wt  \sS,D_\rho)$)  for some $K'>0$.
		By Proposition \ref{p:acent}(d), we obtain the desired estimate.
		\item 
		Let $x,y \in X$ and $u,v \in \sS_n$ be as in the statement of the lemma. We choose $n_0$ be the integer such that $c(x,y) \subset \sS_{n_0}$. Let $\wt u, \wt v \in \sS_{n_0}$ be the vertices such that $u,v$ are descendants of $\wt u, \wt v$ respectively. Let $\wt c \in c(x,y)$. As shown in the proof of Lemma \ref{l:centhyp}, either $\wt u$ (resp. $\wt v$) is equal to $\wt c$ or is a horizontal neighbor of $\wt u$ (resp. $\wt v$). Therefore by \hyperlink{H2}{$\on{(H2)}$}, $\pi(\wt u) \vee \pi(\wt v) \le K_0 \pi(c(x,y))$. We now construct the desired path $\gamma$ from $u$ to $v$ as follows. We join $u$ to $\wt u$ and $\wt v$ to $v$ using the geneology. We can connect $\wt u$ to $\wt v$ using $\wt c$ if necessary. By  \hyperlink{H1}{$\on{(H1)}$} and \hyperlink{H2}{$\on{(H2)}$}, the length $L_\rho(\gamma)$ of $\gam$ is bounded by
		\[
		L_\rho(\gamma) \le \pi(c(x,y)) \left( 1 + \frac{2K_0}{1-\eta_+}\right).
		\]
	\end{enumerate}
\end{proof}

The following proposition provides a bound on visual parameter on $\partial(\sS,D_\rho)$ and relies crucially on \hyperlink{H3}{$\on{(H3)}$}. This construction  of metric is slightly different from that of \cite{Car} and \cite[Theorem 5.1]{Sha}.

\begin{prop}[Visual parameter control] \label{p:visual}
	Let $(X,d)$ and $\sS$ be as in the statement of Theorem \ref{t:H1-4}.
	Let $\rho:\sS \to (0,\infty)$ which satisfies hypotheses  \hyperlink{H1}{$\on{(H1)}$}, \hyperlink{H2}{$\on{(H2)}$}, and \hyperlink{H3}{$\on{(H3)}$}. Then there exists a visual metric $\theta_\rho$ on $\partial(\sS,D_\rho)$ with base point $v_0 \in \sS_0$ with visual parameter $e$, where $D_\rho$ is as defined in \eqref{e:defDrho}. There exists $C_1>1$ such that the metric $\theta_\rho$ satisfies
	\[
 C_1^{-1} \pi(c(x,y)) \le	\theta_\rho(\mathfrak{l}_\infty^{-1}(x),\mathfrak{l}_\infty^{-1}(y)) \le \pi(c(x,y)) \q \mbox{for all $x,y \in X$.}
	\]
	Furthermore, the  map $p:\left(\partial(\sS,D_\rho), \theta_\rho\right) \to (X,d)$ is a power quasisymmetry.
\end{prop}
\begin{proof}
	We define the desired metric  $\theta_\rho$ as
	\[
	\theta_\rho(\mathfrak{l}_\infty^{-1}(x),\mathfrak{l}_\infty^{-1}(y)) = \inf \Biggl\{\sum_{i=0}^{k-1} \pi(c(x_i,x_{i+1})): k \in \bN, x_0=x, x_k=y, x_i \in X \mbox{ for all $i$} 	\Biggr\}.
	\]
	Clearly, $\theta_\rho$ is non-negative, satisfies the triangle inequality and $\theta_\rho(\mathfrak{l}_\infty^{-1}(x),\mathfrak{l}_\infty^{-1}(y)) \le \pi(c(x,y))$ for all $x,y \in X$. It suffices to show 
	\[
	\theta_\rho(\mathfrak{l}_\infty^{-1}(x),\mathfrak{l}_\infty^{-1}(y)) \gtrsim \pi(c(x,y))
	\]
	To this end consider a sequence $x_0,\ldots,x_k$ such that $x_0=x,x_k=y$. Without loss of generality, we may assume that $x_i \neq x_{i+1}$ for all $i=0,\ldots,k-1$. We choose $n$ large enough so that we can apply Lemma \ref{l:Drho}(d) to each pair $x_i,x_{i+1}, i=0,\ldots,k-1$. Choose $v_i \in \sS_n$ such that $x_i \in B_{v_i}$ for all $i=0,\ldots,k$. By concatenating all points obtained by applying Lemma \ref{l:Drho}(d) to each pair $x_i,x_{i+1}$, we obtain a path $\gam \in \Gam_n(x,y)$ such that
	\[
	L_\rho(\gam) \le C \sum_{i=0}^{k-1} \pi(c(x_i,x_{i+1})),
	\]
	where $C$ is the constant from Lemma \ref{l:Drho}(d). Combining the above estimate with \hyperlink{H3}{$\on{(H3)}$}, we obtain
	\be \label{e:dvis}
	(K_1C)^{-1}\pi(c(x,y))\le	\theta_\rho(\mathfrak{l}_\infty^{-1}(x),\mathfrak{l}_\infty^{-1}(y)) \le \pi(c(x,y)),
	\ee
	where $K_1$ is the constant in \hyperlink{H3}{$\on{(H3)}$}.
	
	By Proposition \ref{p:bbs}(c), the map $p: (\partial(\sS,D_1),\theta_1) \to (X,d)$ is a bi-Lipschitz map where $\theta_1$ is a visual metric on $\partial(\sS,D_1)$ with base point $v_0 \in \sS_0$ and visual parameter $a$. Since the identity map $\on{Id}:(\sS,D_\rho) \to (\sS,D_1)$ is a quasi-isometry by  Lemma \ref{l:Drho}(a), the induced boundary map (as defined in Proposition \ref{p:BS}(b)) $\partial \on{Id}: (\partial(\sS,D_\rho),\theta_\rho) \to (\partial(\sS,D_1),\theta_1)$ is a power quasisymmetry by  Proposition \ref{p:BS}(c).  Composing this power quasisymmetry  $\partial \on{Id}: (\partial(\sS,D_\rho),\theta_\rho) \to (\partial(\sS,D_1),\theta_1)$  with the bi-Lipschitz map $p: (\partial(\sS,D_1),\theta_1) \to (X,d)$  yields the desired conclusion that
	$p:\left(\partial(\sS,D_\rho),\theta_\rho\right) \to  (X,d)$ is a power quasisymmetry.
\end{proof}

Consider the metric $\Theta_\rho:X \times X \to [0,\infty)$ defined by
\be \label{e:defTheta} \Theta_\rho(x,y):=\theta_\rho(\mathfrak{l}_\infty^{-1}(x),\mathfrak{l}_\infty^{-1}(y)).
\ee
where $\theta_\rho$ is the visual metric from Proposition \ref{p:visual}. 
The following is an immediate corollary of Proposition \ref{p:visual}.
\begin{cor} \label{c:visual}
	Let $(X,d)$ and $\sS$ be as in the statement of Theorem \ref{t:H1-4}.
	Let $\rho:\sS \to (0,\infty)$ which satisfies hypotheses  \hyperlink{H1}{$\on{(H1)}$}, \hyperlink{H2}{$\on{(H2)}$}, and \hyperlink{H3}{$\on{(H3)}$}. Then the metric $\Theta_\rho$ defined in \eqref{e:defTheta} satisfies $\Theta_\rho \in \sJ_p(X,d)$.
\end{cor}

The following lemma provides a sequence in $\mathfrak{l}_\infty^{-1}(x) \in \partial (\sS,D_2)$ with desirable properties and is useful for approximating balls centered at $x \in X$ in different metrics.
\begin{lem} \label{l:infinitegeodesic}
	Let $(X,d)$ be a compact doubling metric space
	and let $(\sS,D_2)$ denote the corresponding hyperbolic filling with vertical and horizontal parameters $a,\lam$ respectively such that $a \ge \lam \ge 6$.
	For any $x \in X$, there exists a sequence of vertices $v_n \in \sS_n, n \ge 0$ such that $v_n$ is the parent of $v_{n+1}$ and $d(x,\pi_1(v_n)) \le \frac{1}{1-a^{-1}}a^{-n}$ for all $n \in \bN_{\ge 0}$. 
\end{lem}
\begin{proof}
	For each $n$, we choose $w_n \in \sS_n$ be such that $d(w_n,x) < a^{-n}$ (this is possible since $X_n$ is a maximal $a^{-n}$-separated subset).
	Consider the sequence of genealogies $g(w_n)$ for each $n \in \bN_{\ge 0}$. By a diagonal argument the sequence of genealogies  $g(w_n)$ converge along a subsequence to yield the a sequence $v_n \in \sS_n$ such that $v_{n}$ is the parent of $v_{n+1}$ for all $n \in \bN_{\ge 0}$.
	If $v_n \in \sS_n$ is in the genealogy of $w_k, k >n$, we have $d(x,\pi_1(v_n)) <a^{-k} +\sum_{i=n}^k a^{-k}$. Letting $k \to \infty$ along a subsequence yields the desired bound $d(x,\pi_1(v_n)) \le \frac{1}{1-a^{-1}}a^{-n}$.
\end{proof}

The following lemma provides an approximation of balls in $(X,\Theta_\rho)$ using the balls in $(X,d)$. In the following lemma, we use the notation $B_{\Theta}(\cdot,\cdot), B_d(\cdot,\cdot)$ to denote the balls in the metrics $\Theta_\rho, d$ respectively. 
\begin{lem} \label{l:ballcomp}
	Let $(X,d)$ and $\sS$ be as in the statement of Theorem \ref{t:H1-4}.
	Let $\rho:\sS \to (0,\infty)$ which satisfies hypotheses  \hyperlink{H1}{$\on{(H1)}$}, \hyperlink{H2}{$\on{(H2)}$}, and \hyperlink{H3}{$\on{(H3)}$}. Let $\Theta_\rho \in \sJ_p(X,d)$ be as defined in \eqref{e:defTheta} and let $K_0$ denote the constant in \hyperlink{H2}{$\on{(H2)}$}. 
	Let  $L>1$ such that
	\be \label{e:ad1}
	\frac{1}{L}\pi(c(x,y)) \le \Theta_\rho(x,y) \le \pi(c(x,y)) \quad \mbox{for all $x,y \in X$.}
	\ee
	 For any $x \in X$, let   $(v_n)_{n \in \bN_{\ge 0}}$ denote a sequence such that $v_n \in \sS_n$ for all $n \in \bN_{\ge 0}$ as given in Lemma \ref{l:infinitegeodesic}. Then 
	\be \label{e:ad2}
	B_d(\pi_1(v_k), 2 a^{-k}) \subset B_\Theta(x,r),\quad \mbox{whenever $v_k$ satisfies $\pi(v_k) < K_0^{-1} r$,}
	\ee 
	and
	\be  \label{e:ad3}
	B_\Theta(x,r) \subset B_d(\pi_1(v_k), 2(\lam +2) a^{-k}),  
	\ee
	whenever $k \in \bN$ is such that $\pi(v_k) \le K_0 Lr$ and $\pi(v_{k-1}) > K_0 Lr$.
	
\end{lem}

\begin{proof}
First, we show \eqref{e:ad2}. Let $k \in \bN_{\ge 0}$ be such that $\pi(v_k) < K_0^{-1}r$ and let $y \in B_d(\pi_1(v_k),2 a^{-k})$. Let $(z,l) \in c(x,y)$. Since $\{x,y\} \subset  B_d(\pi_1(v_k),2 a^{-k})$, we have $l \le k$. Note that $D_2((z,l),v_l) \le 1$, since 
$d(\pi_1(v_l),x) \le d(\pi_1(v_l),\pi_1(v_k))+d(\pi_1(v_k),x) \le (1-a^{-1})^{-1} a^{-l}+ (1-a^{-1})^{-1} a^{-k} \le {2 (1-a^{-1})^{-1} a^{-l} < \lam a^{-l}}$ and $x \in B(z,2a^{-1}) \subset B(z,\lam a^{-l})$.
This along with \hyperlink{H2}{$\on{(H2)}$} implies that
\[
\pi(c(x,y)) \le K_0 \pi(v_l) \le K_0 \pi(v_k) < r. 
\]
This estimate along with \eqref{e:ad1} implies \eqref{e:ad2}.

Next, we show \eqref{e:ad3}. Let $k \in \bN$ be such that $\pi(v_k) \le K_0 Lr$ and $\pi(v_{k-1}) > K_0 Lr$ and let $y \in B_\Theta(x,r)$. By \eqref{e:ad1}, $\pi(c(x,y)) \le Lr$.
Let $(z,l) \in c(x,y)$. Note that $D_2(v_l,(z,l))\le 1$ which along with \hyperlink{H2}{$\on{(H2)}$} implies that
\[
\pi(v_l) \le K_0 \pi(c(x,y)) \le K_0 L r.
\]
The choice of $k$ implies that $l \ge k$. Hence by Lemma \ref{l:infinitegeodesic}
\begin{align}
	d(\pi_1(v_k),y) &\le d(\pi_1(v_k),\pi_1(v_l)) + d(\pi_1(v_l),z) + d(z,y) \nonumber \\
	& < (1-a^{-1})^{-1}a^{-k} + 2 \lam a^{-l} + 2 a^{-l} \quad (\mbox{since $D_2(v_l,(z,l))\le 1$ and $(z,l) \in c(x,y)$})\nonumber\\
	&\le \left(2 + (1-a^{-1})^{-1} + 2 \lam \right) a^{-k} < 2(2+\lam) a^{-k}\quad \mbox{(since $k \ge l$).}
\end{align}
This completes the proof of \eqref{e:ad3}.
\end{proof}.

\subsection{Construction of  homogeneous measures using weights} \label{ss:hom}

Next, we need to control the Assouad dimension of $(X,\Theta_\rho)$, where $\Theta_\rho \in \sJ_p(X,d)$ is as given in Corollary \ref{c:visual}. To this end, we construct a $p$-homogeneous measure on $(X,\Theta_\rho)$  using hypothesis \hyperlink{H4}{$\on{(H4)}$}. 
  This along with Theorem \ref{t:VK} implies an upper bound on the Assouad dimension   $\dim_{\on{A}}(X,\Theta_\rho) \le p$. To this end, we construct a doubling measure on $X$ using the weight function $\rho: \sS \to  (0,\infty)$. 
  The idea is to construct a measure on $X$ as a limit of discrete measures on $\mathcal{S}_k$ as $k \to \infty$.
  To this end, we introduce the following notions.
  \begin{definition} \label{d:weight}
  	Let $k \in \bN,p \in (0,\infty), C \in (1,\infty)$ and $f_0:\mathcal{S}_k \to (0,\infty), f_1: \mathcal{S}_{k+1} \to (0,\infty)$. Let $\pi: \mathcal{S} \to (0,\infty)$ be a weight.
  	\begin{enumerate}[(a)]
  		\item We say that $f_0:\mathcal{S}_k \to (0,\infty)$ is \emph{$(C,\pi)$-balanced} if $\frac{f_0(u)}{\pi(u)^p} \le C^2 \frac{f_0(v)}{\pi(v)^p}$ for all vertices $u,v \in \mathcal{S}_k$ with $D_2(u,v)=1$. 
  		\item  A function $f_0:\mathcal{S}_k \to (0,\infty)$ is \emph{$(C,\pi)$-unbalanced on   $e=\{u,v\}$} if either  $\frac{f_0(u)}{\pi(u)^p} > C^2 \frac{f_0(v)}{\pi(v)^p}$ or  $\frac{f_0(v)}{\pi(v)^p} > C^{2} \frac{f_0(u)}{\pi(u)^p}$. Similarly, we say that a function $f_0:\mathcal{S}_k \to (0,\infty)$ is \emph{$(C,\pi)$-balanced on   $e=\{u,v\}$} if $ C^{-2} \frac{f_0(v)}{\pi(v)^p} \le \frac{f_0(u)}{\pi(u)^p} \le C^2 \frac{f_0(v)}{\pi(v)^p}$.
  		\item We say that the pair $(f_0,f_1)$ is \emph{$(C,\pi)$-compatible} if for all points $u \in \sS_k$ and $v \in \sS_{k+1}$ such that $u$ is the parent of $v$, we have
  		\be \label{e:mure2}
  		\frac{f_{0}(u)}{\pi(u)^p} \le    \frac{f_{1}(v)}{\pi(v)^p} 
  		\le  C \frac{f_{0}(u)}{\pi(u)^p}.
  		\ee
  	\end{enumerate}
  \end{definition}
  We remark that the notions of balanced and compatibility depend only on the horizontal and vertical edges of $(\mathcal{S},D_2)$ respectively.

Given a horizontal edge $e=\{u,v\}$ in $\sS_k$, we define the \emph{$(C,\pi)$-balancing operator} $B^e_k:(0,\infty)^{\mathcal{S}_k} \to  (0,\infty)^{\mathcal{S}_k}$ as follows.  If $f_0$ is $(C,\pi)$-balanced on $e$, we set $B_k^e(f_0) \equiv f_0$.  
Otherwise, if $\frac{f_0(u)}{\pi(u)^p} > C^2 \frac{f_0(v)}{\pi(v)^p}$
we set
\[
(B^e_k (f_0))(w)= \begin{cases}
	f_0(w) & \mbox{if $w \notin \{u,v\}$,}\\
	f_0(u)- \alpha_1 & \mbox{if $w = u$,}\\
	f_0(v)+ \alpha_1 & \mbox{if $w =v$,}
\end{cases}
\]
where $\alpha_1$ is given by 	\[
\alpha_1 \left( \frac{C^2}{\pi(v)^p}+ \frac{1}{\pi(u)^p} \right) = \frac{f_0(u)}{\pi(u)^p}- C^2 \frac{f_0(v)}{ \pi(v)^p},
\]
so that $\frac{	(B^e_k (f_0))(u)}{ \pi(u)^p} = C^2 \frac{	(B^e_k (f_0))(v)}{ \pi(v)^p}$.  	The case $\frac{f_0(v)}{\pi(v)^p} > C^{2} \frac{f_0(u)}{\pi(u)^p}$ is similar to $\frac{f_0(u)}{\pi(u)^p} > C^2 \frac{f_0(v)}{\pi(v)^p}$.  The terminology is due to the fact that $B_k^e(f_0)$ is $(C,\pi)$-balanced on $e$ for all $f \in (0,\infty)^{\mathcal{S}_k}$.

  We need the following  modification of  a  lemma of Vol'berg and Konyagin 	\cite[Lemma, p. 631]{VK} which plays a key role in the construction of doubling measures.  
\begin{lemma}\label{l:redistribute}
	Let $(X,d)$ and $\sS$ be as in the statement of Theorem \ref{t:H1-4}. Let $\rho:\sS\to (0,\infty)$ be function that satisfies hypotheses \hyperlink{H1}{$\on{(H1)}$}  and  \hyperlink{H4}{$\on{(H4)}$}.
	Let $C \ge \eta_-^{-p}$, where the constants $\eta_-,p$ are as given in the hypotheses \hyperlink{H1}{$\on{(H1)}$} and \hyperlink{H4}{$\on{(H4)}$}.
	Let $k \in \bN_{\ge 0}$, and let $\mu_k$ be a probability mass function on $\sS_k$ such that $\mu_k$ is $(C,\pi)$-balanced.
	Then there exists a 
	probability mass function $\mu_{k+1}$ on $\sS_{k+1}$   such that  the following hold.
	\begin{enumerate}[\rm (1)]
		\item The pair $(\mu_k,\mu_{k+1})$ is $(C,\pi)$-compatible.
		\item The function $\mu_{k+1}$ is $(C,\pi)$-balanced. 
		\item The construction of the measure $\mu_{k+1}$ from the measure 
		$\mu_k$ can be regarded as the transfer of masses from the points of
		$X_k$ to those of $X_{k+1}$, with no mass transferred over a distance greater 
		than $(1+2 \lam a^{-1})a^{-k}$. More precisely, there is a probability measure $\mu_{k,k+1}$ on $X \times X$ which is a coupling of the probability measures $\wt \mu_k:= \sum_{u \in \sS_k} \mu_k(u) \delta_{\pi_1(u)},  \wt \mu_{k+1}:= \sum_{v \in \sS_{k+1}} \mu_{k+1}(v) \delta_{\pi_1(v)}$ such that 
		\[
		\mu_{k,k+1}\left( \{(x_1,x_2) \in X \times X: d(x_1,x_2) \ge (1 + 2\lam a ^{-1}) a^{-k}\}\right) =0,
		\]
		where $\delta_x$ denotes the Dirac  measure at $x\in X$. Here by a coupling we mean the projection maps from $X \times X$ to $X$ on the first and second component pushes forward the measure $\mu_{k,k+1}$ to $\wt \mu_k$ and $\wt\mu_{k+1}$ respectively.
	\end{enumerate}
\end{lemma}

The proof of Lemma \ref{l:redistribute} is done in two steps. First is an `averaging' step where we construct a measure on $\sS_{k+1}$ by distributing the mass $\mu_k(u)$ of every vertex  $u \in \sS_k$  to its children so that the mass received by each child $v$ is proportional to $\pi(v)^p$. At end of this step, we obtain a measure which satisfies the compatibility condition but not necessarily $(C,\pi)$-balanced. In the second `balancing' step, we ensure that the measure is $(C,\pi)$-balanced by a repeated \emph{local transfer} of mass along edges in $\sS_{k+1}$ using the balancing operators $B^e_{k+1}$. By a \emph{local} transfer we mean that the mass is transferred from a vertex to its neighbor. The next two lemmas show useful properties of the balancing operators. The first one shows that the compatibility condition is preserved by balancing operators.
\begin{lem} \label{l:balance}
	Let $(X,d)$ and $\sS$ be as in the statement of Theorem \ref{t:H1-4}. Let $\rho:\sS\to (0,\infty)$ be function that satisfies hypotheses \hyperlink{H1}{$\on{(H1)}$}  and  \hyperlink{H4}{$\on{(H4)}$}. Let  $\mu_k$ be a probability mass function on $\sS_k$ that is $(C,\pi)$-balanced for some $C >1$. Let $f_0:\sS_{k+1} \to (0,1]$ be a probability mass function on $\sS_{k+1}$ such that $(\mu_k,f_0)$ is $(C,\pi)$-compatible. Let $e=\{w_1,w_1'\}$ be an edge in $\sS_{k+1}$ such that $f_0$ is $(C,\pi)$-unbalanced on $e$.
 Then the pair $(\mu_k,B_{k+1}^e(f_0))$ is also $(C,\pi)$-compatible.
\end{lem}
\begin{proof}
  Without loss of generality, we assume that 
		$\frac{f_0(w_1)}{\pi(w_1)^p} > C^2 \frac{f_0(w_1')}{\pi(w_1')^p}$.
Let $v_{1}$ and $v_{1}'$ be parents of $w_{1},w_{1}'$ respectively and let $f_1:=B_{k+1}^e(f_0)$. 
	By construction, we have
	\begin{equation}\label{e:in4p}
		f_{1}(w_{1}) < f_0(w_{1}), \hspace{4mm} f_{1}(w_{1}') > f_0(w_{1}').
	\end{equation}
	Therefore by the $(C,\pi)$-compatibility of $(\mu_k,f_0)$ and \eqref{e:in4p}, we have
	\begin{equation*}
		\frac{f_{1}(w_{1})}{\pi(w_{1})^p} \le   C \frac{\mu_k(v_{1})}{\pi(v_{1})^p},  \hspace{4mm}  \frac{f_{1}(w_{1}')}{\pi(w_{1}')^p} \ge  \frac{\mu_k(v_{1}')}{\pi(v_{1}')^p}.
	\end{equation*}
	Therefore it suffices to verify that
	\begin{equation}\label{e:inAB}
		\frac{f_{1}(w_{1})}{\pi(w_{1})^p} \ge   \frac{\mu_k(v_{1})}{\pi(v_{1})^p},  \hspace{4mm}  \frac{f_{1}(w_{1}')}{\pi(w_{1}')^p} \le C \frac{\mu_k(v_{1}')}{\pi(v_{1}')^p}.
	\end{equation}
	Suppose the first inequality in \eqref{e:inAB} fails to be true, then by construction, \eqref{e:in4p} and  the $(C,\pi)$-compatibility of $(\mu_k,f_0)$, we have
	\begin{equation}
		\frac{\mu_k(v_{1})}{\pi(v_{1})^p} > \frac{f_{1}(w_{1})}{\pi(w_{1})^p}  = C^2 \frac{f_{1}(w_{1}')}{\pi(w_{1}')^p} > C^2  \frac{f_0(w_{1}')}{\pi(w_{1}')^p} \ge C^2 \frac{\mu_k(v_{1}')}{\pi(v_{1}')^p},
	\end{equation}
	which implies $\frac{\mu_k(v_{1})}{\pi(v_{1})^p}  > C^2 \frac{\mu_k(v_{1}')}{\pi(v_{1}')^p}$. However,  Lemma \ref{l:filling}(b) implies that $D_2(v_{1},v'_{1})\le 1$ and therefore the above estimate contradicts the assumption that $\mu_k$ is $(C,\pi)$-balanced. This proves the first inequality in \eqref{e:inAB}. The proof of the second inequality in \eqref{e:inAB} is similar. Indeed, assume to the contrary that  $\frac{f_{1}(w_{1}')}{\pi(w_{1}')^p} >    C \frac{\mu_k(v_{1}')}{\pi(v_{1}')^p}$; then we have
	\begin{equation}
		\frac{\mu_k(v_{1})}{\pi(v_{1})^p} \ge C^{-1} \frac{f_{0}(w_{1})}{\pi(w_{1})^p} > C^{-1}\frac{f_{1}(w_{1})}{\pi(w_{1})^p} = C \frac{f_{1}(w_{1}')}{\pi(w_{1}')^p}> C^2  \frac{\mu_k(v_{1}')}{\pi(v_{1}')^p},
	\end{equation}
	which again implies $\frac{\mu_k(v_{1})}{\pi(v_{1})^p}  > C^2 \frac{\mu_k(v_{1}')}{\pi(v_{1}')^p}$, a contradiction to  the assumption that $\mu_k$ is $(C,\pi)$-balanced.
 In particular $(\mu_k,B_{k+1}^e(f_0))$ is $(C,\pi)$-balanced.
\end{proof}

The next property is that a balancing operator cannot create unbalanced edges. More precisely, we have the following lemma.
\begin{lem} \label{l:balance2}
		Let $(X,d)$ and $\sS$ be as in the statement of Theorem \ref{t:H1-4}. Let $\rho:\sS\to (0,\infty)$ be function that satisfies hypotheses \hyperlink{H1}{$\on{(H1)}$}  and  \hyperlink{H4}{$\on{(H4)}$}. Let  $\mu_k$ be a probability mass function on $\sS_k$ that is $(C,\pi)$-balanced for some $C >1$. Let $f_0:\sS_{k+1} \to (0,1]$ be a probability mass function on $\sS_{k+1}$ such that $(\mu_k,f_0)$ is $(C,\pi)$-compatible. Let $e=\{w_1,w_1'\}$ be an edge in $\sS_{k+1}$ such that $f_0$ is $(C,\pi)$-unbalanced on $e$. If an edge $e'=\{w,w'\}$ on $\sS_{k+1}$ is such that $f_0$ is $(C,\pi)$-balanced on $e'$, then $B_{k+1}^e(f_0)$ is also $(C,\pi)$-balanced on $e'$.
\end{lem}

\begin{proof}
	Without loss of generality, we assume that 
	$\frac{f_0(w_1)}{\pi(w_1)^p} > C^2 \frac{f_0(w_1')}{\pi(w_1')^p}$.
	Let $e'=\{w,w'\}$ be such that $f_0$ is $(C,\pi)$-balanced on $e'$.  
	Let $f_1:= B_{k+1}^e(f_0)$.
 	By our assumption  $\frac{f_0(w_1)}{\pi(w_1)^p} > C^2 \frac{f_0(w_1')}{\pi(w_1')}$, we have  $\set{w,w'} \neq \{w_1,w_1'\}$.
	If $\set{w,w'} \cap \{w_1,w_1'\} = \emptyset$, then   there is nothing to prove since $f_0$ and $f_1$ agree on $\{w,w'\}$.

	The only remaining case   to consider is  if $\set{w,w'} \cap \{w_1,w_1'\}$ contains exactly one element.
	Next, we consider the case $\{w_1,w_1'\} \cap \set{w ,w'}= \set{w_{1}}$ where $w_{1}=w$. Since $f_{0}(w)/\pi(w)^p> C^2 f_{0}(w_{1}')/\pi(w_{1}')^p$,
	\begin{equation}\label{e:in7}
		\frac{f_{1}(w)}{\pi(w)^p}=	\frac{f_{1}(w_{1})}{\pi(w_{1})^p} = C^2  \frac{f_{1}(w_{1}')}{\pi(w_{1})^p},\hspace{4mm} f_{1}(w) < f_{0}(w), \hspace{4mm} f_{1}(w')=f_0(w').
	\end{equation}
	We need to show that 
	\begin{equation} \label{e:in6}
			\frac{f_{1}(w')}{\pi(w')^p} \le C^2 \frac{f_{1}(w)}{\pi(w)^p}, \quad \frac{f_{1}(w)}{\pi(w)^p}	 \le C^2 \frac{f_{1}(w')}{\pi(w')^p}.
	\end{equation}
	Therefore by \eqref{e:in7}, only the first inequality in \eqref{e:in6} can fail for $f_{1}$. Suppose that this happens, that is
	\begin{equation}\label{e:in8}
		\frac{f_{1}(w')}{\pi(w')^p} > C^2 \frac{f_{1}(w)}{\pi(w)^p}.
	\end{equation}
	Let $v',v_{1}'\in \sS_k$ be parents of $w',w_{1}'$ respectively.
	By Lemma \ref{l:balance}, $(\mu_k,f_1)$ is $(C,\pi)$-compatible.
	Then by the $(C,\pi)$-compatibility of $(\mu_k,f_1)$, \eqref{e:in8}, and \eqref{e:in7}, we obtain
	\begin{equation}
		\frac{\mu_k(v')}{\pi(v')^p} \ge C^{-1}\frac{f_{1}(w')}{\pi(w')^p} \stackrel{\eqref{e:in8}}{>} C  \frac{f_{1}(w)}{\pi(w)^p} \stackrel{\eqref{e:in7}}{=} C^3 \frac{f_{1}(w_{1}')}{\pi(w_{1}')^p} \ge C^3 \frac{\mu_k(v_{1}')}{\pi(v_{1}')^p}> C^2 \frac{\mu_k(v_{1}')}{\pi(v_{1}')^p},
	\end{equation}
	which contradicts the assumption that $\mu_k$ is $(C,\pi)$-balanced (since $D_2(v',v_{1}') \le 1$ by Lemma \ref{l:filling}(b) and $\lam \ge 2 + 4 \lam a^{-1}$). 
	The remaining case  $\{w_1,w_1'\} \cap \set{w ,w'}= \set{w_{1}'}$  is analyzed similarly and therefore the assertion that $B_{k+1}^e(f_0)$ is also $(C,\pi)$-balanced on $e'$ is proved.
\end{proof}

The following iterative construction uses Lemmas \ref{l:balance} and \ref{l:balance2} to obtain a balanced and compatible function from a compatible function. 
\begin{lem} \label{l:ind}
		Let $(X,d)$ and $\sS$ be as in the statement of Theorem \ref{t:H1-4}. Let $\rho:\sS\to (0,\infty)$ be function that satisfies hypotheses \hyperlink{H1}{$\on{(H1)}$}  and  \hyperlink{H4}{$\on{(H4)}$}. Let  $\mu_k$ be a probability mass function on $\sS_k$ that is $(C,\pi)$-balanced for some $C >1$. Let $f_0:\sS_{k+1} \to (0,1]$ be a probability mass function on $\sS_{k+1}$ such that $(\mu_k,f_0)$ is $(C,\pi)$-compatible. Let 
		$p_i=\{v_i,v_i'\}, i=1,\ldots,T$ be an enumeration of all edges in $\sS_{k+1}$. 
		We inductively define
		\begin{equation} \label{e:induct}
			f_i := B^{p_i}_{k+1}(f_{i-1}) \in (0,1]^{\sS_{k+1}} \quad \mbox{for all $i=1,\ldots,T$.}
		\end{equation} 
		Then, $f_i, i=0,\ldots,T$ satisfy the following properties.
		\begin{enumerate}[(a)]
			\item Each $f_i$ is a probability mass function such that $(\mu_k,f_i)$ is $(C,\pi)$-compatible for all $i=0,1,\ldots,T$.
			\item The probability mass function $f_T$ is $(C,\pi)$-balanced.
			\item There are no pairs of edges $p_l=\set{w_1,w_2}$, $p_{n}=\set{w_2,w_3}$, $l,n \in \bZ \cap [1,T], l < n$, such that mass is transferred from $w_1$ to $w_2$  in the transition from $f_{l-1}$ to $f_{l}$  and then mass is transferred from $w_2$ to $w_3$  in the transition from $f_{n-1}$ to $f_{n}$.
		\end{enumerate}
\end{lem}
\begin{proof}
	\begin{enumerate}[(a)]
		\item Since the balancing operators preserve the sum, each $f_i$ is a probability mass function. By Lemma \ref{l:balance}, $(\mu_k,f_i)$ is $(C,\pi)$-compatible.
		\item This is an immediate consequence of Lemma \ref{l:balance2}, since $f_T$ is $(C,\pi)$-balanced on every edge in $\sS_{k+1}$.
		\item  Assume the opposite; that is, there are a mass transfer from
		$w_1$ to $w_2$ (in the transition from $f_{l-1}$ to $f_{l}$) followed by a mass transfer from $w_2$ to $w_3$ (in the transition from $f_{n-1}$ to $f_{n}$ with $n>l$), so that
		\begin{equation} \label{e:in9}
			\frac{f_{l}(w_1)}{\pi(w_1)^p} = C^2 \frac{f_{l}(w_2)}{\pi(w_2)^p}, \hspace{4mm} \frac{f_{n-1}(w_2)}{\pi(w_2)^p}> C^2 \frac{f_{n-1}(w_3)}{\pi(w_3)^p}.
		\end{equation}
		By choosing $l$ as the largest number less than $n$ such that   mass is transferred into $w_2$ in the transition from $f_{l-1}$ to $f_{l}$ before the mass transfer from $w_2$ to $w_3$ takes place in the transition from $f_{n-1}$ to $f_{n}$, we may assume that
		\be \label{e:cnd2} 
		f_{l}(w_2)=f_{n-1}(w_2).
		\ee
		If $v_1,v_3$ denote the parents of $w_1,w_3$ respectively, then by Lemma \ref{l:filling}(b)
		\[
		D_2(v_1,v_3) \le 1.
		\]
		Consequently by assumption that $\mu_k$ is $(C,\pi)$-balanced, we have $\mu_k(v_1)/\pi(v_1)^p \le C^2 \mu_k(v_3)/\pi(v_3)^p$. 
		However   \eqref{e:in9}, $(C,\pi)$-compatibility of $(\mu_k,f_l),(\mu_k,f_{n-1})$ along with \eqref{e:cnd2} imply the opposite inequality  $\mu_k(v_1)/\pi(v_1)^p > C^2 \mu_k(v_3)/\pi(v_3)^p$. 
		We have arrived at the desired contradiction and therefore  the property (c) is verified.
	\end{enumerate}
\end{proof}
% \proof  By   \hyperlink{H4}{$\on{(H4)}$}
% \be \label{e:ind0}
% \sum_{w: v \mbox{ is the parent of $w$}} \rho(w)^p \le 1
% \ee
% for all $v \in \sS$. Since $\rho(w) \ge \eta_-$ for all $w$, we obtain that the number of children of each vertex is uniformly bounded above by $S:=\eta_-^{-p}$.
% We choose
% \[
% C=C_1 S,
% \]
% where $C_1:=K_0 \vee \eta_-^{-1}$ is chosen such that (using \hyperlink{H1}{$\on{(H1)}$} and  \hyperlink{H2}{$\on{(H2)}$})
% \be \label{e:ind1}
% \pi(v) \le C_1 \pi(w)
% \quad \mbox{whenever $v,w \in \sS$ satisfies $D_2(v,w)\le 1$.} \ee

Next, we prove  Lemma \ref{l:redistribute} by using the inductive construction in Lemma \ref{l:ind}.

\noindent \textit{Proof of Lemma \ref{l:redistribute}.}
	Let $k \in \bN_{\ge 0}$, and let  $\mu_k$ be a $(C,\pi)$-balanced probability mass function on $\sS_k$.  
	As explained earlier, the transfer of mass is accomplished in two steps.
	In the first `averaging' step, we distribute the mass $\mu_k(v)$ to all its children  such that the mass distributed to each child $w$  is proportional to $\pi(w)^p$ (or equivalently $\rho(w)^p$); that is
	\[
	f_0(w)= \frac{\pi(w)^p}{\sum_{w' \in C(v)} \pi(w')^p} \mu_k(v),
	\]
	for all $v \in \sS_k$ and $w \in C(v)$, where $C(v)$ is as defined in \eqref{e:defchild}.

	By  \hyperlink{H4}{$\on{(H4)}$}, \hyperlink{H1}{$\on{(H1)}$} and the fact that every vertex has at least one child, we obtain 
	\[
	\eta_-^p \pi(v)^p\le \sum_{w' \in C(v)}  \pi(w')^p \le \pi(v)^p.
	\]
	Therefore, we have
	\begin{equation} \label{e:in2}
		\frac{\mu_{k}(v)}{\pi(v)^p} \le    \frac{f_0(w)}{\pi(w)^{p}} \le  \eta_-^{-p}\frac{\mu_{k}(v)}{\pi(v)^p}  \le C \frac{f_0(w)}{\pi(w)^{p}}
	\end{equation}
	for all points $v \in \sS_k$ and $w \in C(v)$. Note that every point $v \in \sS_k$ has at least one child, because we always have $(\pi_1(v),\pi_2(v)+1) \in C(v)$ for any $v \in \sS$. This implies $f_0$ is probability mass function on $\sS_{k+1}$ such that $(\mu_k,f_0)$ is  $(C,\pi)$-compatible as shown in \eqref{e:in2}.
	
	Let $f_T$ denote the probability mass function constructed from $\mu_k$ and $f_0$ as given by Lemma \ref{l:ind}. We claim that $\mu_{k+1}:=f_T$ is the  probability mass function on $\sS_{k+1}$ with the desired properties. Next, we show that $\mu_{k+1}$ satisfies the conditions. 
	\begin{enumerate}[(1)]
		\item 
	 This is an immediate consequence of Lemma \ref{l:ind}(a).
	 \item 
	 This follows from Lemma \ref{l:ind}(b).
	\item
	It remains to verify condition (3).  Since  $d(\pi_1(v),\pi_1(w)) <a^{-k}$ for all $w \in C(v), v \in \sS_k$, there was a mass transfer over a distance of at most $a^{-k}$ while passing from $\mu_k$ to $f_{0}$. Therefore it suffices to verify that while passing from $f_0$ to $f_T=\mu_{k+1}$ there  is a transfer over a distance of at most $2 \lam a^{-k-1}$. Since $d(\pi_1(w),\pi_1(w')) < 2 \lambda a^{-k-1}$ for all points $w,w' \in \sS_{k+1}$ such that $D_2(w,w')=1$, the desired conclusion follows from Lemma \ref{l:ind}(c).
\end{enumerate}
\qed

We construct a doubling measure on $(X,d)$ in Lemma \ref{l:doubling}   using Lemma \ref{l:redistribute}.
\begin{lem}[Construction of doubling measure] \label{l:doubling}
	Let $(X,d)$ be a compact doubling metric space and let $\sS$ denote a hyperbolic filling with vertical and horizontal parameters $a,\lam$ respectively with $a \ge \lam \ge 6$. Let $\rho: \sS \to (0,\infty)$ denote a weight function that satisfies the hypotheses \hyperlink{H1}{$\on{(H1)}$}, \hyperlink{H2}{$\on{(H2)}$}, and \hyperlink{H4}{$\on{(H4)}$}.
	Let $\mu_0$ denote the (unique) probability measure on $\sS_0=\{v_0\}$. Let $\mu_k$ denote the probability measure on $\sS_k$ for all $k \in \bN$ constructed inductively using Lemma \ref{l:redistribute}. Let 
	\[
	\wt \mu_k= \sum_{v \in \sS_k} \mu_k(v) \delta_{\pi_1(v)} \quad \mbox{for all $k \in \bN_{\ge 0}$,}
	\]
	denote a sequence of probability measures on $X$ associated with the above construction. Then any sub-sequential weak limit $\mu$  of $(\wt \mu_k)_{k \in \bN}$ is a doubling measure on $(X,d)$.
\end{lem}
\begin{proof}
	Observe that such a sub-sequential limit $\mu$ exists by Prokhorov's theorem along with the compactness of $(X,d)$.
	
	Since $\diam(X,d)=\frac 1 2$, it suffices to consider $r <1$. 
	For $x \in X$ choose a sequence $\{v_n\}$ as given in Lemma \ref{l:infinitegeodesic}. 
	We obtain two sided bounds on $\mu(B(x,r))$ using $\mu_n(v_n)$ for a suitably chosen value of $n$.
	To describe this let $n \in \bN_{\ge 0}$ denote the largest integer such that $a^{-n} \ge r$. We claim that 
	\begin{equation} \label{e:mub1}
		\mu(B(x,r)) \asymp \mu_n(v_n)
	\end{equation}
	where the constants of comparison are independent of $x \in X$, $r \in (0,1)$. Let us first show the upper bound. If mass from $\mu_n(v), v \in \sS_n$ contributes to $\mu (\overline{B(x,r)})$, then by Lemma \ref{l:redistribute}(3) 
	we have
	\[
	d(\pi_1(v),x) \le r+ \sum_{k=n}^\infty (1+2 \lam a^{-1}) a^{-k}= \left( 1+(1+2 \lam a^{-1}) (1-a^{-1})^{-1} \right) a^{-n}
	\]
	Since $\lam >((1-a^{-1})^{-1}) (1+ 2 \lam a^{-1})$, we have that $x \in B(v,\lam a^{-n}) \cap B(v_n,\lam a^{-n})$ and hence $D_2(v,v_n)\le 1$.  
	Therefore 
	\[
	\mu(B(x,r)) \le \sum_{v \in \sS_n: D_2(v,v_n) \le 1} \mu_n(v).
	\]
	If $D_2(v,v_n) \le 1$ and $v \in \sS_n$, then by Lemma \ref{l:redistribute}(1) and \hyperlink{H2}{$\on{(H2)}$}, we obtain $\mu_n(v_n) \asymp \mu (v)$ for any pair of such vertices. Furthermore since $(X,d)$ satisfies the metric doubling property, the number of neighbors of each vertex is uniformly bounded above \cite[Proposition 4.5]{BBS}. Combining the above estimates yields the upper bound in \eqref{e:mub1}.
	
	For the lower bound, we consider $\mu_{n+2}(v_{n+2})$. By Lemma \ref{l:redistribute}(3) and $d(\pi_1(v_{n+2}), x )<(1+2 \lam a^{-1})a^{-(n+2)}$, we note that the mass from $v_{n+2}$ stays within $B(x, s )$ where \[s= (1+2 \lam a^{-1}) \left( 1 + (1-a^{-1})^{-1} \right) a^{-(n+2) } < r\](since $a^{-n-1}<r$ and $a^{-1}(1+2 \lam a^{-1})(1-a^{-1})^{-1}<1$).
	This implies that $\mu(B(x,r)) \ge \mu_{n+2}(v_{n+2})$. This along with Lemma \ref{l:redistribute}(2) and \hyperlink{H1}{$\on{(H1)}$}, we obtain $\mu_{n+2}(v_{n+2}) \asymp \mu_n(v_n)$. Combining these estimates yields the lower bound for $\mu(B(x,r))$ in \eqref{e:mub1}.
	
	Next, we show that \eqref{e:mub1} implies the desired doubling property.
	For the remainder of the proof we assume $r \in (0,1/2)$. The case $r \ge 1/2$ is similar and easier.
	Let $N \in \bN_{\ge 0}$ denote th largest integer such that $a^{-N} \ge 2r$. This implies $a^{-(N+2)} < 2 a^{-1}r< r$. This implies that $n=N$ or $n=N+1$. Therefore by the same argument as above (using Lemma \ref{l:redistribute}(2) and \hyperlink{H1}{$\on{(H1)}$}), we have $\mu_{n}(v_n) \asymp \mu_{N}(v_N)$. This along with \eqref{e:mub1} shows that $\mu$ is a doubling measure on $(X,d)$.
\end{proof}

Let $\Theta_\rho$ denote the metric defined in \eqref{e:defTheta}.
In the following proposition, we obtain upper bound on the Assouad dimension of $(X,\Theta_\rho)$. We establish this by showing that the measure $\mu$ in Lemma \ref{l:doubling} is $p$-homogeneous in $(X,\Theta_\rho)$. This along with Theorem \ref{t:VK} shows that $\dim_{A}(X,\Theta_\rho) \le p$.
\begin{prop}\label{p:assouad}
	Let $(X,d)$ be a compact doubling metric space and let $\sS$ denote a hyperbolic filling with vertical and horizontal parameters $a,\lam$ respectively such that $a \ge \lam \ge 6$. Let $\rho: \sS \to (0,\infty)$ denote a weight function that satisfies the hypotheses \hyperlink{H1}{$\on{(H1)}$}, \hyperlink{H2}{$\on{(H2)}$},  \hyperlink{H2}{$\on{(H3)}$}, and \hyperlink{H4}{$\on{(H4)}$}. Let $\Theta_\rho$ denote the metric defined in \eqref{e:defTheta} using Proposition \ref{p:visual}. Then the measure $\mu$ defined in Lemma \ref{l:doubling} is $p$-homogeneous in $(X,\Theta_\rho)$, where $p$ is the constant in \hyperlink{H4}{$\on{(H4)}$}. In particular, $\dim_{\on{A}} (X,\Theta_\rho) \le p$.
\end{prop}
\begin{proof}
	For ease of notation, we abbreviate $\Theta_\rho$ by $\Theta$. By \eqref{e:dvis}, there exists $L >1$ such that \eqref{e:ad1} holds.
	Let $\eta_-,\eta_+,K_0$ denote the constants in \hyperlink{H1}{$\on{(H1)}$} and \hyperlink{H2}{$\on{(H2)}$}.

	Next, we show that $\mu$ is $p$-homogeneous in $(X,\Theta)$; that is, there exists $C>1$ such that
	\be \label{e:hom}
	\frac{\mu(B_\Theta(x,r))}{\mu(B_\Theta(x,s))} \le C \left( \frac{r}{s} \right)^p \quad \mbox{ for all $x \in X, 0< s< r$}.
	\ee
	Let $0<s<r$ and $x \in X$. Choose a sequence $(v_n)_{n \in \bN_{\ge 0}}$ such that $v_n \in \sS_n$ for all $n \in \bN_{\ge 0}$ as given in Lemma \ref{l:infinitegeodesic}.  
	Let $k \in \bN_{\ge 0}$ be the smallest non-negative integer such that  $\pi(v_k) < K_0^{-1} s$. By Lemma \ref{l:ballcomp}, Lemma \ref{l:doubling} and \eqref{e:mub1}, we have
	\be \label{e:ad4}
	\mu(B_\Theta(x,s)) \gtrsim \mu_k(v_k).
	\ee
	If $k=0$, then it we have $1 \ge \mu(B_\Theta(x,r)) \ge \mu(B_\Theta(x,s)) \gtrsim 1$ which implies \eqref{e:hom}.
	So it suffices to consider the case $k \ge 1$. The choice of $k$ along with \hyperlink{H1}{$\on{(H1)}$} implies that 
	\be \label{e:ad5} 
	\eta_- K_0^{-1} s \le \eta_- \pi(v_{k-1}) \le \pi(v_k) < K_0^{-1} s.
	\ee
	
	Next, we bound $\mu(B_\Theta(x,r))$ from above. We consider two cases depending on whether or not $r < (K_0L)^{-1}\pi(v_0)$. 
	If $r < (K_0L)^{-1}\pi(v_0)$, there exists 
	$l \in \bN$  such that  $\pi(v_l) \le K_0 Lr$ and $\pi(v_{l-1}) > K_0 Lr$. Hence by \hyperlink{H1}{$\on{(H1)}$}, we have
	\be \label{e:ad7} 
	\eta_- K_0 Lr< \eta_- \pi(v_{l-1}) \le \pi(v_l) \le K_0 Lr.
	\ee
	By Lemma \ref{l:ballcomp}, Lemma \ref{l:doubling}, and \eqref{e:mub1}, we have
	\be \label{e:ad6}
	\mu(B_\Theta(x,r)) \le \mu\left( B_d(\pi_1(v_l), 2(\lam +2) a^{-l}) \right) \lesssim  \mu \left( B_d(\pi_1(v_l),  a^{-l}) \right) \lesssim \mu_l(v_l).
	\ee
	Since $r>s$, we have $l \le k$. Therefore by Lemma \ref{l:redistribute}(2), we have
	\be  \label{e:ad8}
	\frac{\mu_l(v_l)}{\pi(v_l)^p} \le  \frac{\mu_k(v_k)}{\pi(v_k)^p}, \quad \mbox{for any $l \le k$.} 
	\ee
	By \eqref{e:ad4}, \eqref{e:ad6}, \eqref{e:ad5}, \eqref{e:ad7} and \eqref{e:ad8}, we have 
	\[
	\frac{\mu(B_\Theta(x,r))}{\mu(B_\Theta(x,s))} \lesssim \frac{\mu_l(v_l)}{\mu_k(v_k)}  \lesssim \frac{\pi(v_l)^p}{\pi(v_k)^p} \asymp \frac{r^p}{s^p}.
	\]
	This implies \eqref{e:hom} in the case $r < (K_0L)^{-1}\pi(v_0)$.
	
	On the other hand, if $r \ge (K_0L)^{-1}\pi(v_0)$ we use the trivial bound $\mu(B_\Theta(x,r)) \le 1=\mu_0(v_0)$.
	By \eqref{e:ad4},  \eqref{e:ad5},  \eqref{e:ad8} and the bound $ 1 \asymp \pi(v_0) \lesssim r$ , we have 
	\[
	\frac{\mu(B_\Theta(x,r))}{\mu(B_\Theta(x,s))} \lesssim \frac{\mu_0(v_0)}{\mu_k(v_k)}  \lesssim \frac{\pi(v_0)^p}{\pi(v_k)^p} \lesssim \frac{r^p}{s^p}.
	\]
	This completes the proof of \eqref{e:hom}. By Theorem \ref{t:VK}, we obtain the desired bound on Assouad dimension.
\end{proof}

\noindent \emph{Proof of Theorem \ref{t:H1-4}.} This follows immediately from  Corollary  \ref{c:visual} and Proposition \ref{p:assouad}. \qed
\subsection{Upper bound on Assouad dimension using weights}
In this subsection, we prove Theorem \ref{t:simp}.
The proof of the Theorem \ref{t:simp} is very similar to that of \cite[Theorem 1.2]{Car} except for the use of Theorem \ref{t:H1-4} instead of \cite[Theorem 1.1]{Car}. 
For the convenience of the reader, we provide further details since the hypothesis \hyperlink{H4}{$\on{(H4)}$} is different from \hyperlink{H4-old}{$\on{(\wt{H4})}$} of \cite{Car}.
To the reader who is familiar with Carrasco's work, we point out that the estimate in \cite[(2.52)]{Car} implies our version of \hyperlink{H4}{$\on{(H4)}$} for small enough $\eta_0$. The proof of other three hypothesis is  similar. Readers who are familiar with the proof of \cite[Theorem 1.2]{Car} may want to skip the proof of Theorem \ref{t:simp}.

Let $\rho:\sS \to [0,\infty)$ be a function. We define $\rho^*: \sS \to [0,\infty)$ as
\be \label{e:defrhostar}
\rho^*(v) = \min  \{ \rho(w): w \in \sS: \pi_2(w)=\pi_2(v), D_2(v,w) \le 1\} \quad \mbox{for all $v \in \sS$.}
\ee
Similarly, we define $\pi^*: \sS \to [0,\infty)$ as
\be \label{e:defpistar}
\pi^*(v) = \min  \{ \pi(w): w \in \sS: \pi_2(w)=\pi_2(v), D_2(v,w) \le 1\} \quad \mbox{for all $v \in \sS$.}
\ee
If $\gam=(v_1,\ldots,v_N)$ is a horizontal path, we define
\be \label{e:defLh}
L_h(\gam, \rho)= \sum_{j=1}^{N-1} \rho^*(v_j) \wedge \rho^*(v_{j+1}).
\ee
We introduce the following hypothesis on $\rho:\sS \to [0,\infty)$ which serves as a simpler sufficient condition for  \hyperlink{H3}{$\on{(H3)}$}:
		\begin{enumerate}
	\item[ \hypertarget{H3'}{$\on{(H3')}$}] for all $k \ge 1$, for all $v \in \sS_k$ and for all $\gamma \in \Gamma_{k+1}(v)$, it holds $L_h(\gamma,\rho) \ge 1$,
\end{enumerate}
where $L_h(\gam,\rho)$ is as defined in \eqref{e:defLh}.
The hypothesis  \hyperlink{H3'}{$\on{(H3')}$} is simpler to verify than  \hyperlink{H3}{$\on{(H3)}$} because it only involves  curves with horizontal edges.
The following is a version of \cite[Proposition 2.9]{Car} and provides a useful sufficient condition for \hyperlink{H3}{$\on{(H3)}$}.
\begin{prop} \label{p:H3}
	Let $(X,d)$ be a compact doubling metric space.  Let $(\sS,D_2)$ denote the hyperbolic filling with horizontal and vertial parameters $\lam,a$ respectively that satisfy $a \ge \lam \ge 6$. Assume that there exists $p>0$ and a function $\rho: \sS \to (0,\infty)$ which satisfy the hypotheses  \hyperlink{H1}{$\on{(H1)}$}, \hyperlink{H2}{$\on{(H2)}$}, and   \hyperlink{H3'}{$\on{(H3')}$}.
	Then the function $\rho$ also satisfies  \hyperlink{H3}{$\on{(H3)}$}.
\end{prop}
%\noindent \emph{Proof sketch}:
%The proof of Proposition \ref{p:H3} is obtained by making minor modifications to the proof of \cite[Proposition 2.9]{Car} that we outline now. We remark that the assumption (H4) and uniform perfectness in  \cite[Proposition 2.9]{Car}  play no role in its proof.
The proof of Proposition \ref{p:H3} requires several lemmas.
We say that a path $\gam=(v_1,\ldots,v_N)$ is of \emph{level} $k$ (resp. \emph{level at most} $k$) if $\pi_2(v_i)=k$ (resp. $\pi_2(v_i) \le k$) for all $i=1,\ldots,N$.
\begin{lem} (Cf. \cite[Lemma 2.10]{Car} ) \label{l:hor}
	Let $(X,d)$ and $(\sS,D_2)$ be as given in Proposition \ref{p:H3}.
	Let $k \ge 0$ and $v \in \sS_k$. Assume that $\rho$ satisfies \hyperlink{H3'}{$\on{(H3')}$}. Consider a horizontal path $\gam=(v_1,\ldots,v_n)$ of level $k+1$ such that $\pi_1(v_i) \in B(\pi_1(v), 3 a^{-k})$ for all $i=1,\ldots,N$, $\pi_1(v_1) \in B(\pi_1(v),a^{-k})$ and $\pi_1(v_N) \notin B(\pi_1(v), 2 a^{-k})$. Let $w$ denote the parent of $z_1$. Then 
	\[
	\sum_{i=1}^{N-1} \pi^*(v_i)\wedge \pi^*(v_{i+1}) \ge \max \{\pi^*(v),\pi^*(w)\}.
	\]
\end{lem}
\begin{proof}
	First, we show that for all $j=1,\ldots,N$, 
	\be \label{e:hor1}
	\pi^{*}(v_j) \ge  \max \{\pi^*(v),\pi^*(w)\} \min \{ \rho(w_j): w_j \in \sS_{k+1},  D_2(w_j,v_j) \le 1 \}.
	\ee
	Let $\wt w_j \in \sS_{k+1}$ be such that $\pi^*(v_j)= \pi(\wt w_j)$ and $D_2(\wt w_j, v_j) \le 1$. Let $u_j \in \sS_{k}$ be the parent of $w_j$. Then by Lemma \ref{l:filling}(a), 
	\begin{align*}
	d(\pi_1(v),\pi_1(u_j)) &\le d(\pi_1(v), \pi_1(v_j)) +d(\pi_1(v_j),\pi_1(\wt w_j)) + d(\pi_1 (\wt w_j), \pi_1(u_j) )\\
	&< 3 a^{-k} + 2 \lam a^{-k-1}+ a^{-k} =(4+2\lam a^{-1}) a^{-k} < \lam a^{-k},\\
		d(\pi_1(w),\pi_1(v)) & \le 	d(\pi_1(w),\pi_1(v_1))+	d(\pi_1(v_1),\pi_1(v)) < a^{-k} + a^{-k} < \lam a^{-k}
	\end{align*}
The above estimates imply that $D_2(v,u_j) \le 1$ and $D_2(w,u_j) \le 1$.
Therefore $\pi(u_j) \ge \max \{\pi^*(v),\pi^*(w)\}$ and hence
\begin{align*}
\pi^*(v_j) & = \pi(\wt w_j)= \pi(u_j) \rho(\wt w_j) \\
&\ge \max \{\pi^*(v),\pi^*(w)\} \min \{ \rho(w_j): w_j \in \sS_{k+1},  D_2(w_j,v_j) \le 1 \}.
\end{align*}
This completes the proof of \eqref{e:hor1}. Therefore, we have
\begin{align*}
	\sum_{i=1}^{N-1} \pi^*(v_i) \wedge \pi^*(v_{i+1}) &\stackrel{ \eqref{e:hor1}}{\ge}  \max \{\pi^*(v),\pi^*(w)\} 	\sum_{i=1}^{N-1} \rho^*(v_i) \wedge \rho^*(v_{i+1})\\
	& =\max \{\pi^*(v),\pi^*(w)\} L_h(\gam,\rho) \ge  \max \{\pi^*(v),\pi^*(w)\},
\end{align*}
where we use \eqref{e:hor1} in the first line and \hyperlink{H3'}{$\on{(H3')}$} and the second line above.
\end{proof}

We introduce a different notion of length on paths.
For any edge $e=\{u,v\}$ we define
\be  \label{e:deflhat}
\wh \ell_1(e)=\begin{cases}
	\pi^*(u) \wedge \pi^*(v) & \mbox{if $e=\{u,v\}$ is a horizontal edge,}\\
	K_0 \eta_-^{-1} \pi^*(v) & \mbox{if $e=\{u,v\}$ and $u$ is a parent of $v$},
\end{cases}
\ee 
and for a path $\gam= (v_1,\ldots,v_N)$, we define
\be \label{e:deflhatgam}
\wh \ell_1(\gam)= \sum_{i=1}^{N-1} \wh \ell_1(e_i), \quad \mbox{where $e_i=\{v_i,v_{i+1}\}$.}
\ee
If $w \in \sS_{k+1}, u,v \in \sS_k$ such that $D_2(u,v)=D_2(u,w)=1$, then by \hyperlink{H1}{$\on{(H1)}$} and \hyperlink{H2}{$\on{(H2)}$}, we have
\be \label{e:goingdowniscostly}
\wh \ell_1(\{u,v\}) \le \pi^*(u) \le \pi(u) \le  \eta_-^{-1} \pi(w) \le K_0  \eta_-^{-1} \pi^*(w) \le \wh \ell_1(\{u,w\}).
\ee

\begin{lem} \label{l:atmostk}  [Cf. \cite[Lemma 2.11]{Car}]
		Let $(X,d)$ and $(\sS,D_2)$ be as given in Proposition \ref{p:H3}.
	Assume that $\rho: \sS \to [0,\infty)$ satisfies the  hypotheses  \hyperlink{H1}{$\on{(H1)}$}, \hyperlink{H2}{$\on{(H2)}$}, \hyperlink{H3'}{$\on{(H3')}$}. Let $u,v \in \sS_{k+1}$ be such that $d(\pi_1(u), \pi_1(v)) > 4 a^{-k}$. Let $\gam=(v_1,\ldots,v_N)$ be a path of level at most $k+1$ from $v_1=u$ to $v_N=v$. Then there exists a path $\gam'= (u_1,\ldots,u_M)$ of level at most $k$ such that:
	\begin{enumerate}[1.]
		\item $u_1,u_M$ are parents of $v_1$ and $v_N$ respectively, and
		\item $\wh \ell_1(\gam') \le \wh \ell_1 (\gam)$.
	\end{enumerate}
\end{lem}
\begin{proof}
	Let $\gam=(v_1,\ldots,v_N)$ be a path of level at most $k+1$ as given in the statement of the lemma. We decompose $\gam$ into sub-paths of level at most $k$ or level equal to $k+1$. Let $s_1=1$. Define inductively positive integers $s_i, t_i$ as
	\begin{align*}
		t_i&= \min \{ j > s_i: \pi_2(v_j) \le k \mbox{ or } j=N\}, \\
		s_{i+1} &= \min \{j \ge  t_i: \pi_2(v_{j+1}) =k+1\}.
	\end{align*}
We stop when $t_i=N$ for some $i=L$. Note that $\pi_2(v_{s_1})=\pi_2(v_{t_L})=k+1$, and $\pi_2(v_{s_i})= \pi_2(v_{t_j})=k$ for $i \neq 1$ and $j \neq L$.
Since we are trying to bound $\wh \ell_1(\gam)$ from below, we may assume that path $\gam$ has no self-intersections; that is $v_i \neq v_j$ for all $i \neq j$. In particular, $v_{s_i} \neq v_{t_i}$ for all $i$.

For each $i=1,\ldots,L$, let $\gam_i$ denote the sub-path $(v_{s_i},\ldots,v_{t_i-1})$. We will replace each path $\gam_i$ with $\gam_i'$ such that $\wh \ell_1(\gam_i') \le \wh \ell_1(\gam_i)$. 

First, we consider $2 \le i \le L-1$ and postpone the cases $i=1,L$ to the end. Let $2 \le i \le L-1$. We consider two cases.

\noindent \textbf{Case 1:} $\pi_1(v_j) \in B(\pi_1(v_{s_i}), 2 a^{-k})$ for all $j=s_i+1, \ldots,t_i -1$.

In this case, $v_{s_i}, v_{t_i} \in \sS_k$ and are the parents of $v_{s_i+1}, v_{t_i - 1}$ respectively. By Lemma \ref{l:filling}(c), 
$D_2(v_{s_i}, v_{t_i})=1$ and hence we replace $\gamma_i$ with $\gam_i'=(v_{s_i}, v_{t_i})$. From \eqref{e:goingdowniscostly}, we obtain
\[
\wh \ell_1(\gam_i') \le \wh \ell_1 (\{v_{s_i},v_{s_i+1}\}) \le \wh \ell_1(\gam_i).
\]
\noindent \textbf{Case 2:} There exists $j_1 \in \{s_i+1,t_i-1\}$ such that $\pi_1( v_{j_1}) \notin B(\pi_1(v_{s_1}, 2 a^{-k}))$. We assume $j_1$ is the first index with this property.  We denote $j_0=s_i+1, w_{0}=v_{s_i} \in \sS_k$. Suppose $j_l,w_l$ are defined, and if $j_l<t_i-1$, we define
\[
j_{l+1}= \min \{j_l < j <t_i-1:\pi_1 (v_{j}) \notin B(\pi_1(w_l), 2 a^{-k}) \mbox{ or } j=t_i-1\},
\]
and let $w_{l+1} \in \sS_k$ be the parent of $v_{j_l+1} \in \sS_{k+1}$.
Let $L_i$ be such that $j_{L_i}=t_i-1$.

If $l \in \{0,\ldots,L_i-2\}$, we have $\pi_1(v_{j_{l+1}}) \notin B(\pi_1(w_l),2 a^{-k})$.
Since $a > 2 \lam$, we have
$d(\pi_1(w_l), \pi_1(v_{j_l+1})) \le d(\pi_1(w_l), \pi_1(v_{j_l}))+d(\pi_1(v_{j_l}), \pi_1(v_{j_l+1})) < 2 a^{-k} + 2 \lam a^{-k-1}  < 3 a^{-k}$. Therefore by Lemma \ref{l:hor}, Lemma \ref{l:filling}(c), and \eqref{e:goingdowniscostly}, we have
\be \label{e:alm1}
\wh \ell_1 ((w_l,w_{l+1})) \le \pi^*(w_l) \le \wh \ell_1  \left( (v_{j_l},\ldots,v_{j_{l+1}-1}) \right) \quad \mbox{for all $l=0,\ldots,L_i-2$.}
\ee
The above estimate \eqref{e:alm1} is also true for $j=L_i-1$ by combining the above argument and with that of case 1 by considering depending on whether or not $\pi_1( v_{j_1}) \notin B(\pi_1(v_{s_1}, 2 a^{-k}))$.
Hence $\gam_i'=(w_0,\ldots,w_{L_i-1},w_{L_i})$, where $w_0=v_{s_i}, w_{L_i}=v_{t_i}$. By \eqref{e:alm1} along with the above remark, we obtain
\[
\wh \ell_1(\gam'_i) \le \wh \ell_1(\gam_i), \quad \mbox{for $i \in \{2,\ldots,L-1\}$.}
\]
The case $i=1$ is also similar to above. Let $u_1$ be the parent of $v_1$.
Similar to argument above, we consider two cases depending on whether or not $\pi_1(v_j) \in B(\pi_1(u_1),2 a^{-k})$ for all $j=1,\ldots, t_i -1$ as explained in \cite[proof of Lemma 2.11]{Car}. This yields a path $\gam_1'$ from $u_1$ to $v_{t_1}$. The case $i=L$ is exactly same as $i=1$ after reversing the order in which the vertices of $\gam_L$ appear.
By concatenating the paths $\gam_1',\ldots,\gam_L'$, we obtain the path $(u_1,\ldots,u_M)$ with desired properties.
\end{proof}

\begin{lem} (Cf. \cite[Lemma 2.12]{Car}) \label{l:H3}
		Let $(X,d)$ and $(\sS,D_2)$ be as given in Proposition \ref{p:H3}.
		Assume that $\rho: \sS \to [0,\infty)$ satisfies  the  hypotheses  \hyperlink{H1}{$\on{(H1)}$}, \hyperlink{H2}{$\on{(H2)}$}, and \hyperlink{H3'}{$\on{(H3')}$}.  There exists a constant $K_2 \ge 1$ such that the following property: for all $x,y \in X$, there exists $k_0$ depending on $x,y$ such that for all $k \ge k_0$, if $u,v \in \sS_k$ such that $x \in B_u, y \in B_v$, then any path $\gam$ joining $u$ and $v$ satisfies
		\[
		\wh \ell_1(\gam) \ge K_2^{-1} \pi(c(x,y)).
		\]
\end{lem}
\begin{proof}
	Let $u,v \in \sS_k$ be such that $x \in B_u, y\in B_v$.
	Let $m$ be such that $\pi_2(w)=m$ for some (or equivalently, for all) $w \in c(x,y)$. By \eqref{e:boundm}, we have
	$d(x,y) > a^{-m-1}$. For $k \ge m+2$, we have (using $a \ge 12$)
	\be \label{e:lower1}
	d(\pi_1(u),\pi_1(v)) \ge d(x,y)  - d(x,\pi_1(u))- d(y,\pi_1(v)) > a^{-m-1} -2 a^{-m-2} \ge 10 a^{-m-2}.
	\ee
	The idea is to use Lemma  \ref{l:hor} to find a path of level at most $m+2$ whose $\wh \ell_1$ length is larger than $\wh \ell_1(\gam)$. We consider two cases.
	
\noindent \textbf{Case 1}: The path $\gam$ is of level at most $k$, where $k \ge m+2$.
By \eqref{e:lower1}, we can apply Lemma \ref{l:hor}. Set $\gam_k=\gam$. 
Let $u_l, v_l \in \sS_l$  be such that $u,v$ are descendants of $u_l,v_l$ respectively. By Lemma \ref{l:filling}(a), for all $l \ge m+2$, we have
\begin{align*}
d(\pi_1(u_l),\pi_1(v_l)) &\ge 	d(\pi_1(u),\pi_1(v)) - 	d(\pi_1(u),\pi_1(u_l)) - 	d(\pi_1(v),\pi_1(v_l)) \\
&\ge 10 a^{-m-2} -2 \frac{a}{a-1} a^{-l} >6 a^{-m+2}.
\end{align*}
Using the above estimate, and applying Lemma \ref{l:hor} repeatedly we obtain 
path $\gam_{m+2}$ of level at most $m+2$ from $u_{m+2}$ to $v_{m+2}$ such that $\wh \ell_1(\gam) \ge \wh \ell_1 (\gam_{m+2})$. This along with \eqref{e:goingdowniscostly}, \hyperlink{H1}{$\on{(H1)}$}, \hyperlink{H2}{$\on{(H2)}$}, implies
\[
\wh \ell_1(\gam) \ge K_0^{-2} \pi(u_{m+2}) \ge K_0^{-2} \eta_-^{2}  \pi(u_m) \ge K_0^{-3} \eta_-^{2}  \pi(c(x,y)).
\]
In the last estimate, we used $D_2(u_m,w) \le 1$ for any $w \in c(x,y)$ (since $x \in B(\pi_1(u_m), \lam a^{-m}) \cap B(\pi_1(w), \lam a^{-m}) \neq \emptyset$).

\noindent \textbf{Case 2}: $\gam$ is not a path of level at most $k$.
Let $n >k$ be the smallest integer such that $\gam$ is a path of level at most $n$.
Let $k_0 \ge m+2$ be large enough so that
\[
K_0^{-3} \eta_-^{2}\eta_-^m \ge 4 K_0 \eta_-^{-1} \sum_{i=k_0}^\infty \eta_+^i.
\]
Let $\wt u_n, \wt v_n \in \sS_n$ be such that $x \in B_{\wt u_n}, y \in B_{\wt v_n}$ and let $\wt u_k, \wt v_k \in \sS_k$ be the ancestors of  $\wt u_n, \wt v_n$ respectively. By Lemma \ref{l:filling}(a), $D_2(\wt u_n, u_n) \le 1$ and $D_2(\wt v_n, v_n) \le 1$. Let $\gam_u$ denote the path from $\wt u_n$ to $u$ formed by concatenating the genealogy from $\wt u_n$ to $\wt u_k$ and adding an edge from $\wt u_k$ to $u_k$ if necessary. Similarly, let $\gam_v$ denote the path from $v$ to $\wt v_n$ formed in a similar fashion. By concatenating $\gam_u, \gam, \gam_v$ we obtain a path $\wt \gam$ from $\wt v_n$ to $\wt u_n$ whose level is at most $n$. Using the first case, we obtain 
\[
\ell_1 (\wt \gam) \ge K_0^{-3} \eta_-^{2} \pi(c(x,y)) \ge K_0^{-3} \eta_-^{2}\eta_-^m \ge 4 K_0 \eta_-^{-1} \sum_{i=k_0}^\infty \eta_+^i \ge 2 \ell_1(\gam_u) + 2 \ell_1 (\gam_v).
\]
This implies 
\[
\ell_1(\gam) \ge \frac{1}{2}K_0^{-3} \eta_-^{2} \pi(c(x,y))
\]
for any $k \ge k_0$.
\end{proof}

\noindent \emph{Proof of Proposition \ref{p:H3}.}
By \hyperlink{H1}{$\on{(H1)}$},  \hyperlink{H2}{$\on{(H2)}$}, there exists $c>0$ such that 
\[
L_\rho(\gam) \ge c \wh \ell_1(\gam) \quad \mbox{for all paths $\gam$ in $(\sS,D_2)$.}
\]
This estimate along with  Lemma \ref{l:H3} implies \hyperlink{H3}{$\on{(H3)}$}.
\qed

 The statement of the lemma below is slightly different from that of \cite[Lemma 2.13]{Car} and the proof is omitted as it is similar to \cite{Car}.
\begin{lem}(\cite[Lemma 2.13]{Car}) \label{l:VK}
		Let $(X,d)$ and $(\sS,D_2)$ be as given in Theorem \ref{t:simp}.
Suppose we have a function $\pi_0:\sS_k \to (0,\infty)$ such that
\be 
\frac{1}{K} \le  \frac{\pi_0(v)}{\pi_0(w)} \le K \quad \mbox{for all $v,w \in \sS_k$ such that $D_2(v,w) \le 1$},
\ee
where $K \ge 1$ is a constant.
Suppose also that there is function  $\pi_1:\sS_{k+1} \to (0,\infty)$ such that for any $u \in \sS_k$ and for any $v \in \sS_{k+1}$ such that $u$ is the parent of $v$, we have 
\be 
1 \le \frac{\pi_0(u)}{\pi_1(v)} \le K.
\ee
Let $\wh \pi_1: \sS_{k+1} \to [0,\infty)$ be defined as 
\be \label{e:defhatpi} 
\wh \pi_1(w) = 
	\pi_1(w) \vee \left( \frac{1}{K} \max \{ \pi_1(v): v \in \sS_{k+1}, D_2(v,w) \le 1 \} \right).
\ee
 Then
for all $w_1,w_2 \in \sS_{k+1}$ such that $D_2(w_1,w_2) \le 1$, we have
	\be 
	\frac{1}{K} \le  \frac{\wh \pi_1(w_1)}{\wh \pi_1(w_2)} \le K.
	\ee
\end{lem}

\begin{lem}(see \cite[Lemma 2.14]{Car}) \label{l:whtau}
	Let $G=(V,E)$ be a graph whose vertices has a degree bounded by $K$ and let $p>0$. Let $\Gam$ be a family of paths of $G$. Let $\tau: V \to [0,\infty)$ that satisfies
	\[
	\sum_{i=1}^{N-1} \tau(v_i) \ge 1 \quad \mbox{for all paths $\gam=(v_1,\ldots,v_N) \in \Gam$.}
	\]
Let $d_G: V \times V \to [0,\infty)$ denote the combinatorial graph distance metric on $V$. Let $\wh \tau: V \to [0,\infty)$ be defined as 
\[
\wh \tau(v) = 2 \max \{\tau(w): w \in V, d_G(w,v) \le 2 \}.
\]
Then
\[
\sum_{i=1}^{N-1} \wh \tau^*(v_i) \wedge \wh \tau^*(v_{i+1}) \ge 1 \quad \mbox{for all paths $\gam=(v_1,\ldots,v_N) \in \Gam$,} 
\]
where $\wh \tau^*(v) = \min \{\wh \tau(w): d_G(w,v) \le 1\}$, and such that
\[
\sum_{v \in V} \wh \tau(v)^p \le 2^p (K^2 +1)\sum_{v \in V} \tau(v)^p.
\]
\end{lem}
The statement of Lemma \ref{l:whtau} is slightly different from that of \cite[Lemma 2.14]{Car} where the term $K^2 +1$ was replaced by $K^2$. This is because the estimate $\#\{  w \in V : d_G(w,v) \le 2 \} \le K^2$ for all $v \in V$ in \cite{Car} must be replaced by  $\#\{  w \in V : d_G(w,v) \le 2 \} \le K^2+1$. The proof is otherwise identical and is omitted.

\noindent \emph{Proof of Theorem \ref{t:simp}.}
Let $\eta_0 \in (0,1)$ whose value will be determined later.
Since $(X,d)$ is doubling there exists $M_1 \in \bN$, depending only on $a, \lam$ and the doubling constant such that the number of neighbors of each vertex in $(\sS,D_1)$ is bounded by $M_1$, and in particular the number of children of each vertex uniformly bounded above \cite[Proposition 4.5]{BBS}. Set 
\[
\eta_-= \left(\eta_0 M_1^{-1}\right)^{1/p} \in (0,1).
\]
Let $\sigma:\sS \to [0, \infty)$ satisfy \hyperlink{S1}{$\on{(S1)}$} and \hyperlink{S2}{$\on{(S2)}$}. We define $\tau: \sS \to [0,\infty)$ as $\tau := (\sigma^p + \eta^p_-)^{1/p} \ge \eta_-$, which also satisfies \hyperlink{S1}{$\on{(S1)}$}. The function $\tau$ satisfies
\be \label{e:simp1}
\sum_{v \in C(u)} \tau(v)^p \le \sum_{v \in C(u)} \left( \sigma(v)^p +\eta_-^p \right) \le 2 \eta_0 \quad \mbox{for all $u \in \sS$}.
\ee
By \cite[Exercise 10.17]{Hei}, there exists $M_2$, depending only on $\lam$ and doubling constant of $(X,d)$, such that 
\be  \label{e:deg}
\#\{w \in\sS: D_2(v,w)=  1, \pi_2(w)=\pi_2(v)\} \le M_2 \quad \mbox{for all $v \in \sS$.}
\ee
That is, the number of horizontal edges at any vertex is uniformly bounded in $M_2$.
By Lemma \ref{l:whtau}, the function $$\wh \tau (v):= 2 \max \{ \tau(w): w \in \sS_{\pi_2(v)}, D_2(v,w) \le 2 \}$$ satisfies condition \hyperlink{(H3')}{$\on{(H3')}$} and
\begin{align} \label{e:simp2}
\sum_{v \in C(u)} \wh \tau(v)^p &\le 2^p (M_2^2+1)\sum_{\substack{w \in \sS_{\pi_2(u)+1}, \\ D_2(w,v) \le 2}}   \tau(w)^p \q \mbox{(by  Lemma \ref{l:whtau})} \nonumber \\
& \le  2^p (M_2^2+1)\sum_{\substack{\wt u \in \sS_{\pi_2(u)}, \\ D_2(w,v) \le 1}} \sum_{w \in C(\wt u)} \tau(w)^p  \q \mbox{(by Lemma \ref{l:filling}(b))} \nonumber \\
& \le  2^{p+1} (M_2^2+1) \sum_{\substack{\wt u \in \sS_{\pi_2(u)}, \\ D_2(w,v) \le 1}} \eta_0 \q \mbox{(by \hyperlink{S2}{$\on{(S2)}$})} \nonumber \\
& \le  2^{p+1} (M_2^2+1) (M_2+1)  \eta_0 \q \mbox{(by \eqref{e:deg})}.
\end{align}
We construct a function $\rho: \sS \to [0,\infty)$ that satisfies 
\begin{enumerate}[(1)]
	\item $\rho(v) \ge \wh \tau(v)$ for all $v \in \sS$.
	\item $\rho$ satisfies  \hyperlink{H2}{$\on{(H2)}$} with $K_0= \eta_-^{-1}$.
	\item $\rho(v) \le \max \{\wh \tau(w): D_2(w,v) \le 1, \pi_2(w)=\pi_2(v)\}$ for all $v \in \sS$.
%	\item \be
%	\sum_{w \in C(v)} \rho(w)^p \le 2^{p+1} (M_2^2+1)^2 \eta_0, \quad \mbox{for all $v \in \sS$.} 
%	\ee
\end{enumerate}
The idea behind the proof is to inductively construct $\rho$ on $\sS_{k}$ for $k=0,1,\ldots$. Since the conditions (2) and (3) depend only on horizontal edges this inductive construction works well.
We pick $\rho(v_0) = \wh \tau (v_0)$ where $v_0 \in \sS_0$. 
Clearly, this satisfies (1), (2), (3) on $\sS_0$ because $\sS_0$ is a singleton set. Suppose we have constructed $\rho$ on $\cup_{j=0}^{i} \sS_{j}$, we construct $\rho$ on $\sS_{i+1}$ as follows. 
Define $\pi_0: \sS_{i} \to (0,\infty), \pi_1: \sS_{i+1} \to(0,\infty)$ as
\[
\pi_0(u)= \prod_{w \in g(u)} \rho(w), \quad \pi_1(v)= \wh \tau(v) \prod_{w \in g(v), w \neq v} \rho(w) =  \wh \tau(v) \pi_0(\wt v)
\]
for all $u \in \sS_i,v \in \sS_{i+1}$, where $\wt v \in \sS_i$ is the parent of $v \in \sS_{i+1}$. Using the estimate $\wh \tau \ge \eta_-$ along with induction hypothesis, $\pi_0,\pi_i$ satisfy the hypotheses of Lemma \ref{l:VK} with $K=\eta_-^{-1}$. Consider the function $\wh \pi_1: \sS_{i+1} \to (0,\infty)$ defined by \eqref{e:defhatpi} as
\[
\wh \pi_1(w) = 
\pi_1(w) \vee \left( \frac{1}{K} \max \{ \pi_1(v): v \in \sS_{i+1}, D_2(v,w) \le 1 \} \right),
\]
and set $\rho: \sS_{i+1} \to (0,\infty)$ as 
\[
\rho(w)= \frac{\wh \pi_1(w)}{\pi_0(\wt w)} \quad \mbox{for all $w \in \sS_{i+1}$, where $\wt w$ is the parent of $w$.}
\]
Since $\wh \pi_1 \ge \pi_1$ the condition (1) is satisfied. By Lemma \ref{l:VK}, the condition (2) above is also satisfied on $\cup_{j=0}^{i+1} \sS_j$. It only remains to check (3) on $\sS_{i+1}$.
For $v \in \sS_{i+1}$, we have two possibilities for $\wh \pi_1(v)$; either $\wh \pi_1(v)=\pi_1(v)$ or $\wh \pi_1(v)= K^{-1} \pi_1(w)$ for some $w \in \sS_{i+1}$ such that $D_2(v,w)=1$. The first possibility implies that $\rho(v)= \wh \tau(v)$ and hence (3) is satisfied for $v$. The other possibility is that $\wh \pi_1(v)= K^{-1} \pi_1(w)  \ge \pi_1(v)$. In this case, let $\wt v, \wt w \in \sS_i$ denote the parents of $v,w$ respectively. By Lemma \ref{l:filling}(b), $D_2(\wt v, \wt w) \le 1$. Therefore by  condition (2) in the induction hypothesis, we have
\[
\rho(v)=  \frac{\wh \tau(w)\pi_0(\wt w)}{K \pi_0(\wt v)} \le \wh \tau(w),
\]
which concludes the proof of condition (3) above. By induction, there exists a function $\rho: \sS \to (0,\infty)$ which satisfies (1), (2), (3) above. 

Next, we want to show that $\rho$ satisfies the upper bound $\rho \le \eta_+$ in \hyperlink{H1}{$\on{(H1)}$} for some $\eta_+\in(0,1)$ and the hypothesis  \hyperlink{H4}{$\on{(H4)}$} whenever $\eta_0$ is small enough. To this end, consider 
\begin{align*}
\sum_{v \in C(u)} \rho(v)^p &\le \sum_{v \in C(u)} \sum_{ \substack{w \in \sS: D_2(v,w) \le 1,\\ \pi_2(v)=\pi_2(w)}} \wh \tau(w)^p \quad \mbox{(by condition (3))}\\
& \le (M_2+1)  \sum_{ \substack{\wt u \in \sS: D_2(u,\wt u) \le 1, \\ \pi_2(\wt u)=\pi_2(u)}} \sum_{v \in C(\wt u)} \wh \tau(v)^p   \mbox{ (by \eqref{e:deg} and Lemma \ref{l:filling}(c))}\\
& \le 2^{p+1} (M_2^2+1) (M_2+1)^3 \eta_0 \quad \mbox{ (by \eqref{e:deg} and \eqref{e:simp2})}.
\end{align*}
By the above estimate, the choice $\eta_0 \in(0,1)$ such that 
$2^{p+1} (M_2^2+1) (M_2+1)^3 \eta_0 = 2^{-p}$ implies the upper bound $\rho \le \eta_+$ in \hyperlink{H1}{$\on{(H1)}$} for $\eta_+= \frac 1 2 \in (0,1)$ and also  \hyperlink{H4}{$\on{(H4)}$}. Since $\wh \tau$ satisfies \hyperlink{(H3')}{$\on{(H3')}$} and $\rho \ge \wh \tau$, $\rho$ also satisfies \hyperlink{(H3')}{$\on{(H3')}$}. This along with conditions (1), (2) above and Proposition \ref{p:H3} implies that $\rho$  satisfies  \hyperlink{H1}{$\on{(H1)}$},  \hyperlink{H2}{$\on{(H2)}$},  \hyperlink{H3}{$\on{(H3)}$}, and  \hyperlink{H4}{$\on{(H4)}$}. 
The desired conclusion follows from Theorem \ref{t:H1-4}.
\qed
\section{Critical exponent associated to the combinatorial modulus} \label{s:modulus}
Let $G=(V,E)$ be a graph and let $\Gam$ be a family of paths in $G$. 
Consider a function $\rho: V \to [0,\infty)$ and for $\gam \in \Gam$, we define its $\rho$-length as
\[
\ell_\rho(\gam) := \sum_{v \in \gam} \rho (v),
\]
and its $p$-mass  by
\[
M_p(\rho)= \sum_{v \in V} \rho(v)^p.
\]
The $p$-combinatorial modulus\footnote{One could alternatively define the function $\rho$ on edges instead of vertices but for bounded degree graphs this would lead to an equivalent quantity. This follows from an argument of He and Schramm  proof in \cite[Proof of Theorem 8.1]{HS}. Our results could be stated in terms of this alternate definition as well.} of $\Gam$ is defined as
\[
\Mod_p(\Gam,G)= \inf_{\rho \in \Adm(\Gam)} M_p(\rho),  
\]
where $\Adm(\Gam):= \{ \rho: V \to [0,\infty)\hspace{1mm} \vert \hspace{1mm}  \ell_\rho(\gam) \ge 1 \mbox{ for all $\gam \in \Gam$}\}$ denote the set of $\Gam$-admissible functions. If $\Gam= \emptyset$, we set $\Mod_p(\Gam, G)= 0$ by convention.

We recall the definition of  critical exponent of the combinatorial modulus associated to a compact metric space $(X,d)$. The idea behind the following definition is to approximate the compact metric space by a sequence of graphs $G_k$. Then the behavior of the modulus of (discrete) family of curves on $G_k$ which `cross an annulus' as $k \to \infty$ determines a critical exponent.  
\begin{definition} \label{d:critical}
Let $a,\lam,L \in (1,\infty)$ and $p>0$ and let $(X,d)$ be a compact metric space. Let $X_k$ denote a maximal $a^{-k}$ separated subset of $X$ for all $k \ge 0$ and let $\sS_k=\{(x,k): x \in X_k\}$. In this section, we need not assume that $X_k$ is increasing in $k$. Let $\pi_1:\sS_k \to X,\pi_2: \sS_k \to \bN_{\ge 0}$ be the projection maps to the first and second components. 
 For each $k \ge 1$, we define a graph $G_k$ whose vertex set is $\sS_k$ and  there is an edge between distinct vertices $v$ and $w$ if and only if $B(\pi_1(v), \lam a^{-\pi_2(v)}) \cap B(\pi_1(v), \lam a^{-\pi_2(v)}) \neq \emptyset$.   For $v \in \sS$, we define
	\begin{equation} \label{e:defGamkL}
	\Gam_{k,L}(v)=\inf \Biggl\{ \gam=(v_1,v_2,\ldots,v_n) \Biggm|
	\begin{minipage}{210 pt}
		$\gam$ is a path in $G_{\pi_2(v)+k}$ with  $\pi_1(v_1) \in B_v,\\ \pi_1(v_n) \notin B(\pi_1(v),L a^{-\pi_2(v)})$	\end{minipage}
	\Biggr\}.
\end{equation}
Define 
\begin{align} \label{e:defMpk}
	M_{p,k}(L) &= \sup_{v \in \sS} \Mod_p(\Gam_{k,L}(v), G_{\pi_2(v)+k}), \nonumber\\
	M_p(L) &= \liminf_{k \to \infty} M_{p,k}(L).
\end{align}
The \emph{critical exponent of the combinatorial modulus} of $(X,d)$ is defined as 
\be \label{e:defQL}
\on{CE}(X,d)= \inf \{p \in (0,\infty): M_p(L)=0\}.
\ee
\end{definition}
If $\rho \in \Adm(\Gam)$, then $1 \wedge \rho \in \Adm(\Gam)$. This shows that $\Mod_p(\Gam,G)$ is non-increasing in $G$ for any family of paths $\Gam$ and any graph $G$. This shows that the set of $p$ such that $M_p(L)=0$ is an interval.

Strictly speaking $\on{CE}(X,d)$ should be denoted as $\on{CE}(X,d,a,\lam,L, \{X_k: k \ge 0\})$ since it might depend on all these choices of $a, \lam, L$ and $\{X_k: k \ge 0\}$.
It is known that this exponent does not depend on the choice of $L>1$ \cite[Lemma 3.3]{Car}. We will show that it also does not depend on the choices of the $a, \lam \in (1,\infty)$ and $\{X_k: k \ge 0 \}$. To this end, we recall the following lemma. Given a set $Y$, we use the notation $2^Y$ and $\# Y$ to denote the power set of $Y$ and the cardinality of $Y$ respectively.
\begin{lem} \cite[Lemma C.4]{Kig22} \label{l:kig}
 Let $G=(V,E), \wt G=(\wt V,\wt E)$ be two  graphs and let $H:V \to 2^{\wt V}$ be a function so that $\# H(v) < \infty$ for all $v \in V$. Let $\Gam,  \wt \Gam$ be two families of paths in $G,\wt G$ respectively such that for each $\gam \in \Gam$, there exists $\wt \gam \in \wt \Gam$ so that $\wt \gam$ is contained in $\cup_{v \in \gam} H(v)$. Then
 \be
 \Mod_p(\Gam ,G) \le \left( \sup_{v \in V} \# H(v) \right)^p  \sup_{\wt v \in \wt V} \# \{v \in V \mid \wt v \in H(v)\}    \Mod_p(\wt \Gam,\wt G).
 \ee
 \end{lem}
The following proposition shows that the critical exponent for the combinatorial modulus is well defined.
\begin{prop} [critical exponent is well defined] \label{p:wd} Let $a,\wt a,\lam, \wt \lam, L, \wt L \in (1,\infty)$. Let $X_k$ (resp. $\wt X_k$) denote a sequence of maximal $a^{-k}$-separated (resp. $\wt a^{-k}$-separated) subsets of $X$. Let $Q, \wt Q$ denote the corresponding critical exponents be as defined in  \eqref{e:defQL} for these two sets of parameters. Then
	\[
	Q=\wt Q.
	\]
\end{prop}
\begin{proof}
	Let $M_{p,k}(L)$ and $\wt M_{p,k}(\wt L)$ be as defined in \eqref{e:defMpk}. Let $G_k, \wt G_k, k \ge 0$ be the corresponding graphs with vertex sets $\sS_k, \wt \sS_k$ respectively.
%Let $p> Q'(L')$. By \cite[Remark 1, p. 534]{Car} and monotonicity of $M'_{p,k}$ in $p$, we have 
%\be \label{e:wd1}
%\lim_{k' \to \infty} M'_{p,k'}=0.
%\ee
By symmetry, it suffices to show that $Q \le \wt Q$. Or equivalently, it suffices to show that 
$Q \le p$ for any $p> \wt Q$. 
To show this, we need an upper bound on $\Mod_p(\Gam_{k,L}(v), G_{\pi_2(v)+k})$ for $v \in G_n, n \in \bN$.
Let $m \in \bZ$ be the unique integer such that
\be  \label{e:wd0}
2 \wt L \wt a^{-m} \le  \frac{(L-1)}{2}a^{-1} < 2 \wt L \wt a^{-m+1}.
\ee 
For any $n \in \bN$, let $\wt n \in \bN$  be the unique positive integer   such that 
\be  \label{e:wd1}
2 \wt L\wt a^{-\wt n + (1-m)_+} \le \frac{(L-1)}{2}a^{-n} < 2 \wt L \wt a^{-\wt n+1+ (1-m)_+}.
\ee 
For any $\wt k \in \bN$, let $k \in \bZ$ be the unique integer such that 
\be \label{e:wd2}
a^{-k} \le  \frac{(L-1)(\wt \lam -1)}{4\wt L \lam \wt a^{1+(1-m)_+}}\wt a^{-\wt k} <a^{-k+1}. 
\ee
It is evident that there exists $k_0 \in \bN$ such that $\wt k \ge k_0$ implies that $k \ge 1$ (that is $k \in \bN$). For the remainder of the proof we assume $\wt k \ge k_0$.
By \eqref{e:wd1} and \eqref{e:wd2}, we have
\be \label{e:wd3}
 \lam a^{-n-k} < (\wt \lam -1) \wt a^{-\wt n -\wt k}, \q  \frac{\wt a^{-\wt n -\wt k}}{a^{-n-k}} \le  \frac{ \wt a a  \lam }{ \wt \lam -1}.
\ee 
For all $l,\wt l \in \bN$, we define a family of maps $H_{l, \wt l}: \sS_{l} \to \wt \sS_{\wt l}$ such that $H_{l, \wt l}(v)= \{w\}$ where $w \in \wt \sS_{\wt l}$ is such that $d(\pi_1(v),\pi_1(w)) < \wt a^{-\wt l}$ or equivalently, $\pi_1(v) \in B(\pi_1(w), \wt a^{-\wt l})$. Since $\bigcup_{u \in \wt \sS_{\wt l}} B(\pi_1(u), \wt a^{-l})=X$ such a $w \in \wt \sS_{\wt l}$ exists.   
 By \cite[Exercise 10.17]{Hei}, there exists $\beta> \dim_A(X,d)$ and $C_1>1$ such that 
\be \nonumber
\sup_{w \in \wt \sS_{\wt l}} \# \{v \in \sS_l : w \in H(v)\} \le C_1 \left(1 \vee \frac{\wt a^{-\wt l}}{a^{-l}}\right)^\beta.
\ee 
In particular, for any $k,\wt k, n,\wt n \in \bN$ that satisfy \eqref{e:wd2} and \eqref{e:wd3}, we have
\be  \label{e:wd4}
\sup_{w \in \wt \sS_{\wt n+\wt k}} \# \{v \in \sS_{n+k} : w \in H(v)\} \le C_1  \left(1 \vee \frac{\wt a a \lam}{\wt \lam -1 } \right)^\beta.
\ee
Let $\gam=(v_1,\ldots,v_N) \in \Gam_{n,k}(v), v \in \sS_n$ denote an arbitrary path. Note that, $v_i \in \sS_{n+k}$ for all $i=1,\ldots,N$. Consider the sequence $(w_1,\ldots,w_N)$ such that $\{w_i\} = H_{n+k, \wt n + \wt k}(v_i)$ for all $i=1,\ldots,N$. By the first estimate in \eqref{e:wd3}, we have $B(\pi_1(v_i), \lam a^{-n-k}) \subset B(\pi_1(w_i), \wt \lam \wt a^{-\wt n -\wt k})$ for all $i$. This in turn implies for any $i=1,\ldots,N-1$, either $w_i=w_{i+1}$ or $w_i$ and $w_{i+1}$ are neighboring vertices in $\wt G_{\wt n +\wt k}$.
 This implies that for any $\gam \in \Gam_{k,L}(v)$ there exists a path $\wt \gam$ in $\wt G_{\wt n + \wt k}$ from $w_1$ to $w_N$. Therefore by the triangle inequality,
 \begin{equation} \label{e:wd5}
 d(\pi_1(w_1),\pi_1(w_N)) \ge (L-1)a^{-n} - 2 \wt a^{-\wt n -\wt k}   
 \stackrel{\eqref{e:wd1}}{\ge}  (4 \wt L - 2 \wt a^{-\wt k}) \wt a^{-\wt n}  \ge 3 \wt L \wt a^{-\wt n} 
 \end{equation}
for all $\wt k \in \bN$ large enough such that $2 \wt a^{-\wt k} \le L$.
Since $\pi_1(v_1) \in B(\pi_1(v),a^{-n}) \cap B(\pi_1(w_1),\wt a^{-\wt n -\wt k})$, we have $d(\pi_1(v), \pi_1(w_1)) < a^{-n}+\wt a^{-\wt n -\wt k}$. There exists $\wt v \in \wt \sS_{\wt n}$ such that $\pi_1(w_1) \in B(\pi_1(\wt v),  \wt a ^{-\wt n})$. Therefore by \eqref{e:wd5}, the path $\wt \gam \in \wt \Gam_{\wt k, \wt L}(\wt v)$ for any $\wt k$ large enough such that $2 \wt a^{-\wt k} \le L$, where 
\begin{align*}
d(\pi_1(\wt v), \pi_1(v)) &\le d(\pi_1(v), \pi_1(w_1))+ d(\pi_1(\wt v), \pi_1(w_1)) \nonumber \\
&< a^{-n}+\wt a^{-\wt n -\wt k}+ \wt  a^{-\wt n} \le a^{-n}+ 2 \wt a ^{-\wt n} \nonumber\\
&  {\le}\wt a ^{-\wt n} \left(2+ \frac{4\wt L   \wt a^{1+(1-m)_+}}{L-1}\right) \q \mbox{(by \eqref{e:wd1})}.
\end{align*}
Setting $\kappa= \left(2+ \frac{4\wt L   \wt a^{1+(1-m)_+}}{L-1}\right)$, we conclude that for all large enough $\wt k \in \bN, n \in \bN, v \in \sS_n, \gamma \in \Gamma_{k,L}(v)$, there exists $\wt \gam \in \wt \Gam_{\wt k, \wt L}$ such that $\wt v \in \wt \sS_{\wt n}$ with $d(\pi_1(\wt v), \pi_1(v))< \kappa \wt a^{-\wt n}$ and $\wt \gam$ is contained in $\bigcup_{u \in \gamma} H_{n+k,\wt n + \wt k}(u)$, where $\wt n,k$ are as given by \eqref{e:wd1} and \eqref{e:wd2} respectively. By \cite[Exercise 10.17]{Hei}, there exists $C_2>1$ such that 
\be \label{e:wd6}
\sup_{v \in \sS_n} \#\{\wt v \in \wt \sS_{\wt n}: d(\pi_1(\wt v), \pi_1(v))< \kappa \wt a^{-\wt n}\} \le C_2 \kappa^\beta.
\ee
Now combining the above with Lemma \ref{l:kig}, \eqref{e:wd4} and \eqref{e:wd6}, we obtain
\begin{align*}
	\Mod(\Gam_{k,L}(v))& \le  C_1  \left(1 \vee \frac{\wt a a \lam}{\wt \lam -1 } \right)^\beta \Mod\left( \bigcup_{ \substack{ \wt v \in \wt \sS_{\wt n}, \\ d(\pi_1(\wt v), \pi_1(v))< \kappa \wt a^{-\wt n} }} \wt \Gam_{\wt k,\wt L}(\wt v) \right) \\
	& \le  C_1  \left(1 \vee \frac{\wt a a \lam}{\wt \lam -1 } \right)^\beta \sum_{  \substack{ \wt v \in \wt \sS_{\wt n}, \\ d(\pi_1(\wt v), \pi_1(v))< \kappa \wt a^{-\wt n} } } \Mod\left(\wt \Gam_{\wt k,\wt L}(\wt v) \right) \\
	& \le  C_1  \left(1 \vee \frac{\wt a a \lam}{\wt \lam -1 } \right)^\beta   C_2 \kappa^\beta \wt M_{p, \wt k}(\wt L)
\end{align*}
for all $n \in \bN, v \in \sS_n$ and for all $\wt k \in \bN$ large enough.
This implies that 
\[
M_{p,k}(L) \le C_1  \left(1 \vee \frac{\wt a a \lam}{\wt \lam -1 } \right)^\beta   C_2 \kappa^\beta \wt M_{p, \wt k}(\wt L)
\]
 for all $\wt k \in \bN$ large enough and for all $k \in \bN$ defined by \eqref{e:wd2}. This immediately implies the desired inequality $Q \le p$ for any $p > \wt Q$.
\end{proof}

The following `reverse volume doubling estimate' is known if the metric space is uniformly perfect \cite[Exercise 13.1]{Hei}. Since our metric space is not necessarily uniformly perfect, the following lemma provides a substitute for uniform perfectness at sufficiently many scales.
\begin{lem} \label{l:rvd}
	Let $\mu$ be a  doubling measure on a metric space  $(X,d)$ such that $\mu(B(x,2r)) \le C_D \mu(B(x,r))$ for all $x \in X,r>0$. Let $x_0,\ldots,x_N$ be a set of points such that $d(x_i,x_{i+1}) <r/4$ for all $i=0,\ldots,N-1$ where $d(x_0,x_N) > R >r$. Then there exists $c, \alpha >0$ depending only on $C_D$ such that 
	$\mu(B(x_0,R)) \ge c (R/r)^{\alpha} \mu(B(x_0,r))$.
\end{lem}
\begin{proof}
	If $s \in [r, R/2]$, then by triangle inequality there exists $x_j$ such that 
 \[
\frac{5}4s \le d(x_0,x_j) \le \frac{7}{4}s.
 \]
 By the doubling property, for such $x_j$, we have
 \[
\mu( B(x_j, s/4) ) \ge C_D^{-3} \mu(B(x_j,4s))  \ge C_D^{-3}\mu(B(x_0,s))
 \] 
 Therefore 
 \[
 \mu(B(x_0,2s)) \ge \mu(B(x_0,s)) + \mu( B(x_j, s/4) )  \ge (1+C_D^{-3}) \mu (B(x_0,s))
 \]
 for all such $s \in [r,R/2]$.
 Let $k$ be the largest  integer such that $2^{k} r \le R$.
 By iterating the above estimate
 \[
 \mu(B(x_0,R)) \ge \mu(B(x_0,2^{k-1} r)) \ge (1+C_D^{-3})^{k-1} \mu(B(x_0,r)) 
 \ge c \left(\frac{R}{r}\right)\mu(B(x_0,r))
 \]
 where $\alpha = \log (1+C_D^{-3})/ \log 2$ and $c= (1+C_D^{-3})^{-1}$.
\end{proof}

\subsection{Proof of Theorem \ref{t:main}}
In this section, we complete the proof of Theorem \ref{t:main}. That is, we show that  for any  compact, doubling, metric space $(X,d)$, we have
\[
\dim_{\on{CA}}(X,d)=\on{CE}(X,d)= \inf \{\dim_{\on{A}}(X,\theta): \theta \in \sJ_p(X,d)\}.
\]
The following lemma is useful to obtain upper bounds on the critical exponent  $\on{CE}(X,d)$.
\begin{lem} \label{l:admissible}
	Let $(X,d)$ be a compact metric space and $a \ge \lambda \ge 6$. Let $\theta \in \sJ(X,d)$ and $\mu$ be a $q$-homogeneous measure on $(X,\theta)$. Let $\mathcal{S}$ denote a hyperbolic filling with parameters $a, \lambda$. For all $v \in \mathcal{S}, k \in \bN$, we define
$\rho_v: \sS_{\pi_2(v)+k} \to [0,\infty)$ as
\be \label{e:defrho}
\rho_v(w) = \begin{cases}
	\left( \frac{\mu(B_w)}{\mu(B_v)}\right)^{1/q} & \mbox{if $B_w \cap  B_d(\pi_1(v), (L +1)a^{-\pi_2(v)})  \neq \emptyset$, $w \in  \gam$ for some $\gam \in \Gam_{k,L}(v)$,}\\
	0 &  \mbox{otherwise.}
\end{cases}
\ee
Then there exists $c>0, k_0 \in \bN$ depending only on $d,\theta,\mu, a, L$ so that
\be \label{e:adm}
\sum_{w \in \gam} \rho_v(w) \ge c \quad \mbox{for all $\gam \in \Gam_{k,L}(v), v \in \sS$}
\ee
for all   $k \ge k_0$.
\end{lem}
\begin{proof}
	 Let $\gam \in \Gam_{k,L}(v)$. To show \eqref{e:adm}, by choosing a sub-path if necessary, we may assume that $\gam=(v_1,\ldots,v_N)$ and $\pi_1(v_1) \in B_v$, $\pi_{1}(v_j)\notin B(\pi_1(v), a^{-\pi_2(v)})=B_v$ for all $j=2,\ldots,N-1$, $\pi_{1}(v_i)\in B(\pi_1(v), L a^{-\pi_2(v)})$ for all $i=1,\ldots,N-1$ and $v_N \notin  \pi_{1}(v_i)\in B_d(\pi_1(v), L a^{-\pi_2(v)})$. 
	Since 
	\be \label{e:ca1}
	d(\pi_1(v_i),\pi_1(v_{i+1})) < 2 \lam a^{-k-\pi_2(v)} \le a^{-k-\pi_2(v)+1} \le a^{-\pi_2(v)} \quad \mbox{for all $i=1,\ldots,N-1$,}
	\ee
	we have   $v_N \in B(v, (L+1) a^{-\pi_{2}(v)})$.
	In particular, 
	\be \label{e:ca2}
	z_i \in B_d(\pi_1(v_i), 2 \lam a^{-\pi_2(v_i)}) \setminus B_d(\pi_1(v_i), a^{-\pi_2(v_i)}) \neq \emptyset  \quad \mbox{for all $i=1,\ldots,N$, }
	\ee
	where $z_i=\pi_1(v_{i+1})$ for $i=1,\ldots,N-1$ and $\pi_1 (v_{i-1})$ for   $i=N$.
	Let $\eta:[0,\infty) \to [0,\infty)$ be a distortion function such that 
	the identity map $\on{Id}:(X,d) \to (X,\theta)$ is an $\eta$-quasisymmetry. This along with the choice of $z_i$ above this implies that
	\be \label{e:ca3}
	B_\theta(\pi_1(v_i), c_1 \theta(\pi_1(v_i), z_i)) \subset B_d(\pi_1(v_i), a^{-\pi_2(v_i)})=B_{v_i} \subset B_\theta(\pi_1(v_i), \eta(1) \theta(\pi_1(v_i), z_i))  \quad 
	\ee
	for all $i=1,\ldots,N$, where $c_1= \left[ \eta \left(  2\lam \right) \right]^{-1}$.
	Since $d( \pi_1(v),\pi_1(v_i)) < (L+1) a^{-\pi_2(v)}$ for all $i=1,\ldots,N$, we have
	$B_v=B_d(\pi_1(v), a^{-\pi_2(v)}) \subset B_d(\pi_1(v_i), (L+2) a^{-\pi_2(v)})$. 
	Choosing $w_i \in \{\pi_1(v), \pi_1(v_N)\}$ such that $2(L+1)a^{-\pi_2(v)} >d(\pi_1(v_i),w_i) \ge d(\pi_1(v), \pi_1(v_N))/2 \ge L a^{-\pi_2(v)}/2$, we have
	\be  \label{e:ca5}
	B_v \subset  B_d(\pi_1(v_i), (L+2) a^{-\pi_2(v)}) \subset B_\theta(\pi_1(v_i), C_2 \theta(\pi_1(v_i),w_i))  
	\ee
	for all $i=1,\ldots,N$, where $C_2=\eta(2(L+2)L^{-1})$. Furthermore
	\be
	B_\theta (\pi_1(v_i),\theta(\pi_1(v_i),w_i)) \subset B_d(\pi_1(v_i), \eta(1) d(\pi_1(v_i),w_i)) \subset B_d(\pi_1(v),C_3 a^{-\pi_2(v)}),
	\ee
	where $C_3=(L+1)(2 \eta(1)+1)$.
	
	Since $\mu$ is $q$-homogeneous on $(X,\theta)$ and $\theta \in \sJ(X,d)$, $\mu$ is a doubling measure on $(X,d)$. Therefore
	\be \label{e:ca4}
	\mu(B_v) \gtrsim \mu\left(   B_d(\pi_1(v),C_3 a^{-\pi_2(v)}) \right) \ge  \mu(B_\theta (\pi_1(v_i),\theta(\pi_1(v_i),w_i)) ).
	\ee 
	Since $d(\pi_1(v_i),\pi_1(v_{i+1}))/d(\pi_1(v_i), w_i)  \le 4L^{-1} \lambda a^{-k}$   
	\be \label{e:ca6}
	\theta(\pi_1(v_i),\pi_1(v_{i+1})) \le \theta(\pi_1(v_i), w_i)    \eta  \left( 4L^{-1} \lam a^{-k}\right)  \q \mbox{for all $i=1,\ldots,N-1$.}
	\ee 
	Pick $k_0 \in \bN$ large enough so that 
	\be  \label{e:ca7}
	c_1    \eta  \left( 4L^{-1} \lam a^{-k}\right) \le 1.
	\ee
	Let $k \ge k_0$. Since $\mu$ is a $q$-homogeneous  measure on $(X,\theta)$ and $\on{Id}:(X,d) \to (X,\theta)$ is an $\eta$-quasisymmetry, we have
	\begin{align*}
		\sum_{i=1}^N \rho_v(v_i) &= 	\sum_{i=1}^N \left( \frac{\mu(B_{v_i})}{\mu(B_v)}\right)^{1/q} \\
		&\gtrsim \sum_{i=1}^{N-1} \left( \frac{\mu(B_\theta(\pi_1(v_i), c_1 \theta(\pi_1(v_i), \pi_1(v_{i+1})))}{ \mu\left( B_\theta(\pi_1(v_i),   \theta(\pi_1(v_i),w_i)) \right)}\right)^{1/q} \q \mbox{(by \eqref{e:ca3} and \eqref{e:ca5})}\\
		& \gtrsim \sum_{i=1}^{N-1} \frac{\theta(\pi_1(v_i), \pi_1(v_{i+1}))}{\theta(\pi_1(v_i), w_i)  } \q \mbox{(by  \eqref{e:ca7} and $q$-homogeneity of $\mu$ in $(X,\theta)$)}\\
		&\gtrsim \sum_{i=1}^{N-1} \frac{\theta(\pi_1(v_i), \pi_1(v_{i+1}))}{\theta(\pi_1(v), \pi(v_N))  }  \q \mbox{(since  $d(\pi_1(v_i),w_i) \lesssim d(\pi_1(v), \pi_1(v_N))$)}\\
		&\gtrsim  \frac{\theta(\pi_1(v_1), \pi(v_N))  }{\theta(\pi_1(v), \pi(v_N))  } \gtrsim 1 \\
		& \quad \mbox{(by triangle inequality and $d(\pi_1(v_1), \pi(v_N)) \gtrsim d(\pi_1(v), \pi(v_N))$).}
	\end{align*}
	This completes the proof of \eqref{e:adm}, where $c>0$ depends only on $\eta$, $q$-homogeneity constants of $\mu$ and $\lam, a, L$.
\end{proof}

\noindent \emph{Proof of Theorem \ref{t:main}.}
 By Proposition \ref{p:wd}, it suffices to consider the critical exponent $a \ge \lam \ge 6$, where the maximal $a^{-n}$ separated subsets are increasing (similar to the definition of hyperbolic filling).

The inequality $\dim_{\on{CA}}(X,d) \le \on{CE}(X,d)$ follows from the same argument as the proof of $\dim_{\on{ARC}}(X,d) \le \on{CE}(X,d)$ in \cite[Theorem 1.3]{Car} where the use of \cite[Theorem 1.2]{Car} is replaced with Theorem \ref{t:simp}. This yields the inequality 
\be \label{e:ex}
\dim_{\on{CA}}(X,d) \le \inf\{\dim_{\on{A}}(X,\theta) : \theta \in \sJ_p(X,d)\} \le\on{CE}(X,d).
\ee

So it suffices to show $\on{CE}(X,d) \le \dim_{\on{CA}}(X,d)$.
Let $p>\dim_{\on{CA}}(X,d)$. We consider a hyperbolic filling $\sS$ of $(X,d)$ with parameters $a \ge \lam \ge 6$. Then by Theorem \ref{t:VK}, for any $q \in (\dim_{\on{CA}}(X,d), p)$, there exists $\theta \in \sJ(X,d)$ and a $q$-homogeneous measure $\mu$ on $(X,\theta)$. 

Next, we show the following estimate: there exists $C>0, k_0 \in \bN$ such that
\be \label{e:bndmod}
\Mod_p \left( \Gam_{k,L}(v)\right) \le C a^{-k \alpha(p-q)} \q \mbox{ for any $v \in \sS, k \in \bN$ with $k \ge k_0$.}
\ee
 Let $\rho_v: \sS_{\pi_2(v)+k} \to [0,\infty)$ be as given in Lemma \ref{l:admissible}. By Lemma \ref{l:admissible}, it suffices to estimate $ \sum_{w} \rho_v(w)^p$.
 
To this end, we first obtain an upper bound on $\rho_v(w)$. For any $v \in \sS$, let  $w \in \sS_{\pi_2(v)+k}$ such that $w \in \gam$ for some $(v_1,\ldots,v_N)= \gam \in \Gam_{k,L}(v)$. 
Note that $d(\pi_1(w),\pi_1(v_1)) \vee d(\pi_1(w),\pi_1(v_N)) \ge \frac 1 2 d(\pi_1(v_1),\pi_1(v_N)) \ge (L-1)a^{-\pi_2(v)}/2$. Therefore there exists a sequence $w=x_0,\ldots,x_M$ so that $d(x_0,x_M) \ge (L-1)a^{-\pi_2(v)}$. Choose $k_1 \in \bN$ such that $(L-1) a^{-k_1}< \frac {1}{4}$.
By the volume doubling property of $\mu$ on $(X,d)$, we have
$\mu(B_v) \gtrsim \mu( B(\pi_1(v), 2L a^{-\pi_2(v)})) \gtrsim \mu( B(\pi_1(w), (L-1) a^{-\pi_2(v)/2}))$.
For all $k \ge k_1$, by Lemma \ref{l:rvd}, we have
\be  \label{e:suprho}
\frac{\mu(B_{w})}{\mu(B_v)} \lesssim \frac{\mu(B_{w})}{\mu( B(\pi_1(w), (L-1) a^{-\pi_2(v)}/2))} \le C a^{-\alpha k}  \quad \mbox{for all $v \in \sS, w\in \gam, \gam \in \Gam_{k,L}(v)$,}
\ee
where $C,\alpha$ only depends on $\lam, a,L,$ and the doubling constant of $\mu$ in $(X,d)$.
 Since $\mu$ is a doubling measure, we have
 \begin{align} \label{e:ca8}
 \sum_{w \in \sS_{\pi_2(v)+k}} \rho_v(w)^q &\le \sum_{\substack{w \in \sS_{\pi_2(v)+k}, \\ \pi_1(w) \in B(\pi_1(v), (L+1) a^{-\pi_2(v)})} }   \frac{\mu(B_w)}{ \mu \left( B_d(\pi_1(v), (L+2) a^{-\pi_2(v)})\right)}  \nonumber \\
 &\lesssim  \sum_{\substack{w \in \sS_{\pi_2(v)+k},\\ \pi_1(w) \in B(\pi_1(v), (L+1) a^{-\pi_2(v)})} } \frac{\mu(B_d(\pi_1(w), a^{-\pi_2(w)}/2))}{ \mu \left( B_d(\pi_1(v), (L+2) a^{-\pi_2(v)}) \right)} \nonumber \\
 & \quad \mbox{(since $\mu$ is doubling)} \nonumber \\
 & \lesssim 1   \q \mbox{ (since $B_d(\pi_1(w), a^{-\pi_2(w)}/2))$ pairwise disjoint).}
\end{align} 		
By  Lemma \ref{l:admissible}, there exists $c >0, k_0 \in \bN$ such that $c^{-1} \rho_v \in \Adm(\Gam_{k,L}(v))$ for al $k \ge k_0$. Hence by \eqref{e:suprho} and \eqref{e:ca8}, we have \[ \Mod_p(\Gam_{k,L}(v)) \lesssim \sum_{w} \rho_v(w)^p  \le \left(\sup_{w} \rho_v(w) \right)^{p-q} \sum_{w} \rho_v(w)^q \lesssim a^{-k \alpha (p-q)} \quad \mbox{for all $k \ge k_0$.} \]
This concludes the proof of \eqref{e:bndmod} and hence we obtain
$M_{p,k}(L) \lesssim a^{-k \alpha (p-q)}$ for all  $k \ge k_0$. This shows $M_p(L)=0$ and hence $\on{CE}(X,d) \le p$ for all  $p>\dim_{\on{CA}}(X,d)$. This along with \eqref{e:ex} concludes the proof. 
\qed

One might wonder if the assumption $\dim_{\on{A}}(X,d)<\infty$ (or equivalently, $(X,d)$ is a doubling metric space) in Theorem \ref{t:main} is necessary. To this end, we present the following example.
\begin{example}
Let $X$ denote the set of all sequences $(x_i)_{i \in \bN}$ such that $x_i \in \{1,2,\ldots,i\}$; that is, $X= \prod_{i=1}^\infty \{1,\ldots,i\}$. We define a metric on $X$ by setting
\[
d((x_i)_{i \in \bN},(y_i)_{i \in \bN} ) = \begin{cases}
	0 & \mbox{if $x_i=y_i$ for all $i \in \bN$,}\\
	2^{-j} & \mbox{if $j= \min\{k: x_k \neq y_k\} <\infty$.}
\end{cases}
\]
It is easy to see that $(X,d)$ is a compact, ultrametric space.
Since every open ball of radius $2^{-k}$ has $k+1$ distinct points with mutual distance of at least $2^{-k-1}$ for each $k \in \bN$, we have \be \label{e:ex1} \dim_{\on{A}}(X,d)=\infty. \ee
On the other hand, consider $a, \lam, L >1, X_k, \sS_k$ as given in Definition \ref{d:critical}. For any $v \in \sS_n$ and any $\gamma=(v_1,\ldots,v_N) \in \Gamma_{k,L}(v)$, since $(X,d)$ is an ultrametric space, we have
\[
(L-1)a^{-n} \le d(\pi_1(v_1),\pi_1(v_N)) \le \max_{1 \le i \le N-1} d(\pi_1(v_i),\pi_1(v_{i+1})) < 2 \lambda a^{-n-k}.
\]
Therefore for all $k$ large enough so that $a^{-k} \le \frac{L-1}{2 \lam}$, we have $\Gam_{k,L}(v) = \emptyset$ and hence $M_{p,k}(L)=0$ for all $p>0$ and $k$ large. This implies $\on{CE}(X,d)=0$. Therefore by \eqref{e:ex1}, we have
\[
 \dim_{\on{CA}}(X,d)=\infty \neq 0=\on{CE}(X,d). 
\]
Therefore the assumption that $(X,d)$ is a doubling metric space is necessary in Theorem \ref{t:main}.
\end{example}
We conclude with some questions about the closely related notion of conformal (Hausdorff) dimension. Recall that the conformal (Hausdorff) dimension $\dim_{\on{CH}}(X,d)$ is defined as
\[
\dim_{\on{CH}}(X,d) = \inf \{  \dim_{\on{H}}(X,\theta): \theta \in \sJ(X,d)\},
\]
where $\dim_{\on{H}}(X,\theta)$ denotes the Hausdorff dimension of $(X,\theta)$. Does the equality  
$$\dim_{\on{CH}}(X,d) = \inf \{  \dim_{\on{H}}(X,\theta): \theta \in \sJ_p(X,d)\}$$ always hold? Theorem \ref{t:main} shows a similar result for the Ahlfors regular conformal dimension.
It is also interesting to know for which metric spaces does $\dim_{\on{CH}}(X,d)= \dim_{\on{CA}}(X,d)$ hold? It is easy to see that $\dim_{\on{CH}}(X,d) \le \dim_{\on{CA}}(X,d)$. One might expect that for `self-similar sets' like the standard Sierpinski carpet $\dim_{\on{CH}}(X,d)= \dim_{\on{CA}}(X,d)$ holds. This seems to be a difficult problem since `self-similarity' is not a quasisymmetry invariant. It is not known whether the equality $\dim_{\on{CH}}(X,d)= \dim_{\on{CA}}(X,d)$  holds even for the standard Sierpinski carpet.
\subsection*{Acknowledgments} 
This material is based upon work supported in part by the National Science Foundation
under Grant No. DMS-1928930 while the author participated in a program hosted
by the Mathematical Sciences Research Institute (MSRI) in Berkeley, California, during
the Spring 2022 semester.

The author benefited from lectures by Mario Bonk and Pekka Pankka on Carrasco Piaggio's work \cite{Car} at MSRI and  discussions after those lectures. I thank Nageswari Shanmugalingam
for informing me about the construction of hyperbolic filling in \cite{BBS,Sha} and K\^ohei Sasaya for comments on an earlier draft of this work. I thank Sylvester Eriksson-Bique and Nageswari Shanmugalingam for pointing out some errors in earlier draft. I am grateful to the anonymous referee for a careful reading and many suggestions that improved the quality of exposition.

\noindent Department of Mathematics, University of British Columbia,
Vancouver, BC V6T 1Z2, Canada. \\
mathav@math.ubc.ca

\end{document}